\documentstyle{article}
\title{State spaces of \jbst s\thanks{Both authors are supported by NSF grant DMS-0101153}}
\author{Matthew Neal\\
Denison University\\
 Granville, OH 43023\\
 nealm@denison.edu \and Bernard Russo\\
University of California\\
Irvine, CA 92697-3875\\
brusso@uci.edu}
\newcommand{\mgt}{minimal geometric tripotent}
\newcommand{\ie}{{\em i.e.}}

\newcommand{\pair}[2]{\langle{#1},{#2}\rangle}

\newcommand{\gtt}{geometric tripotent}
\newcommand{\qed}{\hfill \mbox{$\Box$}}
\newcommand{\ii}{\mbox{$\cal I$}}
\newcommand{\mm}{\mbox{$\cal M$}}

\newcommand{\lc}{$\mbox{sp}$}
\newtheorem{theorem}{Theorem}[section]
\newtheorem{proposition}[theorem]{Proposition}
\newtheorem{corollary}[theorem]{Corollary}
\newtheorem{lemma}[theorem]{Lemma}
\newtheorem{remark}[theorem]{Remark}

\newtheorem{definition}[theorem]{Definition}

\newcommand{\RR}{{\bf R}}
\newcommand{\dd}{{\cal D}}
\newcommand{\fss}{facially symmetric space}
\newcommand{\uu}{\mbox{$\cal U$}}
\newcommand{\ff}{\mbox{$\cal F$}}

\newcommand{\gt}{\mbox{$\cal GT$}}
\renewcommand{\sf}{\mbox{$\cal SF$}}

\newcommand{\CC}{{\bf C}}
\newcommand{\NN}{{\bf N}}
\renewcommand{\t}{\tilde}

\newcommand{\csa}{$C^*$-algebra}
\newcommand{\jbst}{$JB^*$-triple}
\newcommand{\jbwst}{$JBW^*$-triple}

\begin{document}
\maketitle
\begin{abstract}
An atomic decomposition is proved for Banach spaces which satisfy
some affine geometric axioms compatible with notions from the
quantum mechanical measuring process. This is then applied to
yield, under appropriate assumptions,  geometric
characterizations, up to isometry, of the unit ball of the dual
space of a \jbst, and  up to complete isometry, of one-sided
ideals in \csa s. 
\end{abstract}

\tableofcontents

\begin{center}
{\bf Introduction}
\end{center}

 The Jordan algebra of self-adjoint elements of a \csa\ $A$
has long been used as a model for the bounded observables of a
quantum mechanical system, and the states of $A$ as a model for
the states of the system. The state space of this Jordan Banach
algebra is the same as the state space of the \csa\  $A$ and is a
weak$^*$-compact convex  subset of the dual of $A$. With the
development of the structure theory of \csa s, and the
representation theory of Jordan Banach algebras, the problem arose
of determining which compact convex sets in locally convex spaces
are affinely isomorphic to  such a state space. In the context of
ordered Banach spaces, such a characterization has been given for
Jordan algebras in the pioneering paper by Alfsen and Shultz,
\cite{AS}.

After the publication of \cite{AS}, and the corresponding result
for \csa s \cite{AHS}, there began in the 1980s a development of
the theory of \jbst s which paralleled in many respects the
functional analytic aspects of the theory of operator algebras.
\jbst s, which are characterized by holomorphic properties of
their unit ball, form a large class of Banach spaces supporting a
ternary algebraic structure which includes \csa s, Hilbert spaces,
and spaces of rectangular matrices, to name a few examples.
In particular, most of the axioms used by Alfsen and Shultz were
shown to have non-ordered analogs in the context of \jbst s (see
\cite{FriRus85bis}). By the end of the decade, a framework was
proposed by Friedman and Russo in \cite{FriRus89} in which to
study the analog of the Alfsen-Shultz result for \jbst s. A
characterization of those convex sets which occur as the unit ball
of the predual of an {\it irreducible} \jbwst\ was given in
\cite{FriRus93} (see Theorem~\ref{theorem:8.2}). Since \jbst s
have only a local order, the result characterizes the whole unit
ball, which becomes the ``state space'' in this non-ordered
setting.

Guided by the approach of Alfsen and Shultz in the binary context,
it was natural to expect that to prove a geometric
characterization of predual unit balls of {\it global} (that is,
not irreducible) \jbwst s would require a decomposition of the
space into atomic and non-atomic summands and a version of
spectral duality. These goals have remained elusive in the
framework of the axioms used in \cite{FriRus93}. In the present
paper, by introducing the very natural axiom asserting the
existence of a Jordan decomposition in the real linear span of
every norm-exposed face, we are able to prove the atomic
decomposition. In addition, by imposing a spectral axiom every bit
as justified as the one in the Alfsen-Shultz theory, we are able
to give a geometric characterization of the unit ball of the dual
of a  \jbst. These results give positive answers to Problems 1,2
and 3 in \cite{FriRus93}. Moreover, when combined with the recent
characterization of ternary rings of operators (TROs) in terms of
its linear matricial norm structure \cite{NeaRuspj} (see
Theorem~\ref{theorem:5.5}), we obtain a facial operator space
characterization of TROs and one-sided ideals in \csa s, which
responds to a question of D. Blecher.

The main results of this paper are Theorems~\ref{atdec},
\ref{face}, and \ref{C}, which we state here.

\noindent {\bf Theorem~\ref{atdec}.} {\it  Let $Z$ be a neutral,
locally base normed, strongly \fss\ satisfying the pure state
properties and JP. Then $Z=Z_a\oplus^{\ell^1}N$, where $Z_a$ and
$N$ are strongly \fss s satisfying the same properties as $Z$, $N$
has no extreme points in its unit ball, and $Z_a$ is the norm
closed complex span of the extreme points of its unit ball.
Furthermore, by Proposition~\ref{theorem:3.7}, $(Z_a)^*$ is
isometric to an atomic \jbwst.
 }

\medskip

\noindent {\bf Theorem~\ref{face}.} {\it A Banach space $X$ is
isometric to a JB*-triple if and only if $X^*$ is an L-embedded,
locally base normed, strongly spectral, strongly facially
symmetric space which satisfies the pure state properties and JP.
}

\medskip

In \cite{NeaRuspj}, it was proved that an operator space $A$ is
completely isometric to a TRO if and only if $M_n(A)$ is isometric
to a \jbst\ for every $n\ge 2$. Combining this fact with
Theorem~\ref{face} gives a facial operator space characterization
of TRO's. Since a one-sided ideal is a TRO, Theorem~\ref{C} then
gives an operator space characterization of one-sided ideals in
\csa s.

\medskip

\noindent {\bf Theorem~\ref{C}.} {\it Let $A$ be a TRO. Then $A$
is completely isometric to a left ideal in a C*-algebra if and
only if there exists a  convex set $C=\{x_\lambda:\lambda\in
\Lambda\}\subset A_1$  such that the collection of faces
$$F_\lambda:=F_{\left[\begin{array}{c}0\\
x_\lambda/\|x_\lambda\|
\end{array}\right]}\subset M_{2,1}(A)^*,$$ form a directed set with
respect to containment, $F:=\sup_{\lambda}F_\lambda$ exists, and
\begin{description}
\item[(a)] The set $\{\left[\begin{array}{c}0\\ x_{\lambda}
\end{array}\right]:\lambda\in\Lambda\}$ separates the points of $F$;
\item[(b)] $F^\perp=0$ {\rm (}that is, the partial isometry $V\in
(M_{2,1}(A))^{**}$ with $F=F_V$ is maximal{\rm )};
\item[(c)] $\pair{F}{\left[\begin{array}{c}0\\ x_\lambda
\end{array}\right]}\ge 0$ for all $\lambda\in\Lambda$;
\item[(d)] $S_F^*\left(\left[\begin{array}{c}0\\ x_\lambda
\end{array}\right]\right)=\left[\begin{array}{c}0\\ x_\lambda
\end{array}\right]$ for all $\lambda\in\Lambda$.
\end{description}
}

 This paper is organized as follows. In section~\ref{sec:prel} we
recall the background on facially symmetric spaces and on \jbst s
and TRO's. Section~\ref{sec:atomic} is devoted to a proof of the
atomic decomposition. The first subsection contains a result for
some contractive projections on \fss s and the second subsection
introduces and studies the Jordan decomposition property. The
third subsection gives a   geometric characterization of spin
factors (Proposition~\ref{theorem:3.6}), a variation of the main
result of \cite{FriRus92}. The main result of
section~\ref{sec:atomic}, the atomic decomposition
(Theorem~\ref{atdec}), is proved in the fourth subsection.

The main applications occur in section~\ref{sec:ideal}. After
giving a result, interesting in their own right, on contractive
projections on Banach spaces in the first subsection
(Proposition~\ref{theorem:2.4}), the second subsection then uses
all of the machinery developed up to there to give a geometric
characterization of Cartan factors
(Proposition~\ref{theorem:3.7}), a variation of the main result of
\cite{FriRus93}. The spectral duality axiom is introduced in the
next subsection and used together with the atomic decomposition to
give a geometric characterization of the dual ball of a  \jbst\
(Theorem~\ref{face}). The final subsection applies the latter to
give an operator space characterization of one-sided ideals in
\csa s (Theorem~\ref{C}).

\section{Preliminaries}\label{sec:prel}

Facially symmetric spaces (see subsection~\ref{sect:fss}) were
introduced in \cite{FriRus86} and studied in \cite{FriRus89} and
\cite{FriRus92}. In \cite{FriRus93}, the complete structure of
{\em atomic} facially symmetric spaces was determined, solving a
problem posed in \cite{FriRus86}. It was shown, more precisely,
that
 an irreducible, neutral, strongly facially symmetric space is linearly
isometric to the predual of one of the Cartan factors of types 1
to 6, provided that it satisfies some natural and physically
significant axioms, four in number, which are known to hold in the
preduals of all \jbwst s. As in the study of state spaces of
Jordan algebras (see \cite{AS} and the books
\cite{AlfShu01},\cite{AlfShu02}), we shall refer to these axioms
as the pure state properties. Since we can regard the entire unit
ball of the dual of a JB*-triple as the ``state space'' of a
physical system, cf. \cite[Introduction]{FriRus86}, we have given
a geometric characterization of such state spaces.

The project of classifying \fss s was started in \cite{FriRus92},
where, using two of the pure state properties, denoted by STP and
FE, geometric characterizations of complex Hilbert spaces and
complex spin factors were given. The former is precisely a rank 1
\jbwst\ and a special case of a Cartan factor of type 1, and the
latter is the Cartan factor of type 4 and a special case of a
\jbwst\ of rank 2. (For a description of all of the Cartan
factors, see subsection~\ref{sect:jbst}.) The explicit structure
of a spin factor naturally embedded in a facially symmetric space
was then used in \cite{FriRus93} to construct abstract generating
sets and complete the classification in the atomic case.

\subsection{Facially symmetric spaces}\label{sect:fss}
Let $Z$ be a complex normed space.  Elements $f,g\in Z$ are {\em
orthogonal}, notation $f\perp g$, if
$\|f+g\|=\|f-g\|=\|f\|+\|g\|$. A {\em norm exposed face} of the
unit ball $Z_1$ of $Z$ is a non-empty set (necessarily $\neq Z_1$)
of the form $F_x=\{f\in Z_1:f(x)=1\}$, where $x\in Z^*,\|x\|=1$.
Recall that a {\em face} $G$ of a convex set $K$ is a non-empty
convex subset of $K$ such that if $g\in G$ and $h,k\in K$ satisfy
$g=\lambda h+(1-\lambda)k$ for some $\lambda\in (0,1)$, then
$h,k\in G$. In particular, an extreme point of $K$ is a face of
$K$. We denote the set of extreme points of $K$ by ext$\, K$. An
element $u\in Z^*$ is called a {\em projective unit} if $\|u\|=1$
and $\langle u,F_{u}^{\perp}\rangle=0$. Here, for any subset $S$,
$S^\perp$ denotes the set of all elements orthogonal to each
element of $S$.
 \ff\ and \uu\ denote the collections of norm exposed faces of $Z_1$
and projective units in $Z^*$, respectively.

Motivated by measuring processes in quantum mechanics, we defined
a {\em symmetric face} to be a norm exposed face $F$ in $Z_1$ with
the following property: there is a linear isometry $S_F$ of $Z$
onto $Z$, with $S_{F}^2=I$ (we call such maps {\em symmetries}),
such that the fixed point set of $S_F$ is
$(\overline{\mbox{sp}}F)\oplus F^{\perp}$ (topological direct
sum). A complex normed space $Z$ is said to be {\em weakly
facially symmetric} (WFS) if every norm exposed face in $Z_1$ is
symmetric. For each symmetric face $F$ we defined contractive
projections $P_k(F),\ k=0,1,2$ on $Z$ as follows.  First
$P_1(F)=(I-S_F)/2$ is the projection on the $-1$ eigenspace of
$S_F$. Next we define $P_2(F)$ and $P_0(F)$ as the projections of
$Z$ onto $\overline{\mbox{sp}}F$ and $F^{\perp}$ respectively, so
that $P_2(F)+P_0(F)=(I+S_F)/2$. A {\em geometric tripotent} is a
projective unit $u\in\uu$ with the property that $F:=F_u$ is a
symmetric face and $S_{F}^*u=u$ for some choice of symmetry $S_F$
corresponding to $F$. The projections $P_k(F_u)$ are called {\em
geometric Peirce projections}.

\gt\ and \sf\ denote the collections of geometric tripotents and
symmetric faces respectively, and the map $\gt\ni u\mapsto
F_u\in\sf$ is a bijection \cite[Proposition 1.6]{FriRus89}. For
each geometric tripotent $u$ in the dual of a WFS space $Z$, we
shall denote the geometric Peirce projections by $P_k(u)=P_k(F_u),
k=0,1,2$. Also we let $U:=Z^*,Z_k(u)=Z_k(F_u):=P_k(u)Z$ and
$U_k(u)=U_k(F_u):=P_k(u)^*(U)$, so that we have the {\em geometric
Peirce decompositions} $Z=Z_2(u)+Z_1(u)+Z_0(u)$ and
$U=U_2(u)+U_1(u)+U_0(u)$. A symmetry corresponding to the
symmetric face $F_u$ will sometimes be denoted by $S_u$.   Two
geometric tripotents $u_{1}$ and $u_{2}$ are {\em orthogonal} if
$u_{1} \in U_0(u_{2})$ (which implies $u_{2} \in U_0(u_{1})$) and
{\em colinear} if $u_{1} \in U_1(u_{2})$ and $u_{2} \in
U_1(u_{1})$. More generally, elements $a$ and $b$ of $U$ are {\em
orthogonal} if one of them belongs to $U_2(u)$ and the other to
$U_0(u)$ for some geometric tripotent $u$. Two geometric
tripotents $u$ and $v$ are said to be {\em compatible} if their
associated geometric Peirce projections commute, \ie,
$[P_k(u),P_j(v)]=0$ for $k,j\in\{0,1,2\}$. By \cite[Theorem
3.3]{FriRus89}, this is the case if $u\in U_k(v)$ for some
$k=0,1,2$. For each $G\in\ff$, $v_G$ denotes the unique geometric
tripotent with $F_{v_G}=G$.

A contractive projection $Q$ on a normed space $X$ is said to be
{\em neutral} if for each $\xi\in X,\ \|Q\xi\|=\|\xi\|$ implies
$Q\xi=\xi$. A normed space $Z$ is {\em neutral} if for every
symmetric face $F$, the projection $P_2(F)$ corresponding to some
choice of symmetry $S_F$, is neutral.

A WFS space $Z$ is {\em strongly facially symmetric} (SFS) if for
every norm exposed face $F$ in $Z_1$ and every $y\in Z^*$ with
$\|y\|=1$ and $F\subset F_y$, we have $S_F^*y=y$, where $S_F$
denotes a symmetry associated with $F$.

The principal examples of neutral strongly facially symmetric
spaces are preduals of $JBW^*$-triples, in particular, the
preduals of von Neumann algebras, see \cite{Pac}. In these cases,
as shown in \cite{Pac},
 geometric
tripotents correspond to tripotents in a $JBW^*$-triple and to
partial isometries in a von Neumann algebra. Moreover, because of
the validity of the Jordan decomposition for hermitian functionals
on JB*-algebras, $\mbox{sp}_\CC F$ is automatically norm closed
(cf. Lemma~\ref{lem:yfbr}).

In a neutral strongly facially symmetric space $Z$, every non-zero
element has a {\em polar decomposition} \cite[Theorem
4.3]{FriRus89}: for $0\neq f\in Z$ there exists a unique geometric
tripotent $v=v(f)=v_f$ with $f(v)=\|f\|$ and $\langle v,
\{f\}^{\perp}\rangle=0$. Let $\cal M$ denote the collection of
minimal geometric tripotents of $U$,
 {\em i.e.}, ${\cal M}=\{v\in\gt:U_2(v)\mbox{ is one
 dimensional}\}$. If $Z$ is a neutral strongly SFS space
 satisfying PE, then the map $f\mapsto v(f)$ is a bijection of
 ext$\, Z_1$ and $\cal M$ (\cite[Prop. 2.4]{FriRus92}).

A partial ordering can be defined on the set of \gtt s as follows:
if $u,v\in\gt$, then $u\le v$ if $F_u\subset F_v$, or equivalently
(\cite[Lemma 4.2]{FriRus89}), $P_2(u)^*v=u$ or $v-u$ is either
zero or a geometric tripotent orthogonal to $u$.  Let $\cal I$
denote the collection of
 indecomposable geometric tripotents of $U$,
 {\em i.e.}, ${\cal I}=\{v\in\gt:u\in\gt,\ u\leq v\Rightarrow u=v\}$.
 In general, ${\cal M}\subset {\cal I}$, and under certain
conditions, (Proposition\ref{prop:38}(a) below and \cite[Prop.
2.9]{FriRus92}), $\cal M$ coincides with $\cal I$.

We now recall the definitions of the pure state properties and
other axioms.

\begin{definition}\label{defn:2.8}{\rm
 Let $f$ and $g$ be extreme points of the unit ball of
 a neutral SFS space $Z$. The {\em transition probability} of
 $f$ and $g$ is the number
 \[
 \langle f|g\rangle:=f(v(g)).
 \]
 A neutral SFS space $Z$ is said to satisfy {\em ``symmetry of
transition
 probabilities''} STP if for every pair of extreme points
 $f,g\in$  ext$\, Z_1$, we have
 \[
 \overline{\langle f|g\rangle}=\langle g|f\rangle.
 \]
 }\end{definition}

In order to guarantee a sufficient number of extreme points, the
following definition was made in \cite{FriRus92} and assumed in
\cite{FriRus93}. For the present paper, this definition is too
strong and will be abandoned. It will turn out that the property
(b) of Proposition~\ref{prop:1.1} will be available to us and
suffice for our purposes.

\begin{definition}\label{defn:2.5}{\rm
 A normed space $Z$ is said to be {\em atomic} if every symmetric face
 of $Z_1$ has an extreme point.
 }\end{definition}

\begin{definition}\label{defn:2.14}{\rm
 A neutral SFS space $Z$ is said to satisfy property
 FE if every norm closed face of $Z_1$ different from $Z_1$
 is a norm exposed face. We use the terminology PE for the special case of this
 that every extreme point of $Z_1$ is norm exposed.
 }\end{definition}

The following consequence of atomicity will be more useful to us
in this paper.

 \begin{proposition}[\cite{FriRus92},Proposition 2.7]\label{prop:1.1}
 If $Z$ is an atomic SFS space satisfying PE, then
 \begin{description}
 \item[(a)] $U=\overline{\mbox{sp}}\, {\cal M}$ (weak$^*$-closure),
where $\cal M$ is the set of minimal geometric tripotents.
 \item[(b)] $Z_1=\overline{\mbox{co}}\, \mbox{{\em ext}}\, Z_1$ (norm
closure).
 \end{description}
 \end{proposition}

\begin{definition}\label{defn:erp}{\rm
A neutral SFS space $Z$ is said to satisfy the {\em ``extreme rays
property''} ERP if for every $u\in\gt$ and every $f\in$ ext$\,
Z_1$, it follows that $P_2(u)f$ is a scalar multiple of some
element in  ext$\, Z_1$. We also say that $P_2(u)$ preserves
extreme rays.}
\end{definition}

\begin{definition}\label{defn:JP}{\rm
A WFS space $Z$ satisfies JP if for any pair $u,v$ of orthogonal
geometric tripotents, we have
\begin{equation}\label{eq:4.1}
S_uS_v=S_{u+v},
\end{equation}
where for any geometric tripotent $w$, $S_w$ is the symmetry
associated with the symmetric face $F_w$.}
\end{definition}

The property JP was defined and needed in \cite{FriRus93} only for
\mgt s $u$ and $v$. The more restricted definition given here is
needed only in Proposition~\ref{prop:38}(b), where ironically, the
involved geometric tripotents turn out to be minimal. (The
assumption of JP is used in subsection~\ref{subsect:2.1} only for
\mgt s.) As in Remark 4.2 of \cite{FriRus93},with identical
proofs, JP implies the following important {\em joint Peirce}
rules for orthogonal geometric tripotents $u$ and $v$:

\begin{eqnarray*}
Z_{2}(u+v) &=&  Z_{2}(u)+Z_{2}(v)+Z_{1}(u) \cap Z_{1}(v),\\
Z_{1}(u+v) &=&  Z_{1}(u) \cap Z_{0}(v)+Z_{1}(v) \cap Z_{0}(u),\\
Z_{0}(u+v) &=&  Z_{0}(u) \cap Z_{0}(v).\\
\end{eqnarray*}

Definitions~\ref{defn:2.8},\ref{defn:2.14} and \ref{defn:erp} are
analogs of
 physically meaningful axioms in \cite{AS}. 
In the Hilbert space model for quantum mechanics, property JP for
minimal geometric tripotents is interpreted as follows.  Choose
$\xi\otimes\xi$ to be the state exposed by a yes/no question $v$
and $\eta\otimes\eta$ to be the state exposed by another $u$, and
complete $\xi,\eta$ to an orthonormal basis.   For any state
vector $\zeta$ expressed in this basis, the symmetry $S_u$ (resp.
$S_v$) changes the sign of the coefficient of $\xi$ (resp. $\eta$)
and $S_{u+v}$ changes the sign of both coefficients.

We need the concept of $L$-embeddedness for the proofs of
Proposition~\ref{theorem:3.7} and Theorem~\ref{face}. This is
defined as follows. A linear projection  $P$ on a Banach space $X$
is called an {\it $L$-projection} if $\|x\|=\|Px\|+\|(I-P)x\|$ for
every $x\in X$. The range of an $L$-projection is called an {\it
$L$-summand}. The space $X$ is said to be an {\it L-embedded
  space} if it is an
$L$-summand in its second dual. These concepts are studied
extensively in \cite[Chapter IV]{HWW93}. The predual of a \jbst\
is an example of an $L$-embedded space (\cite{BarTim86}) and every
$L$-embedded space is weakly sequentially complete (\cite[Theorem
2.2,page 169]{HWW93}.

 The following is the main result of \cite{FriRus93}. We have
added the assumption of L-embeddedness, which seems to have been
overlooked in \cite{FriRus93}. This omission was discovered in the
process of proving Proposition~\ref{theorem:3.7}. More precisely,
our Proposition~\ref{theorem:2.4} is needed in the proofs of
\cite[Lemmas 5.5 and 6.6]{FriRus93}. In addition, our
Proposition~\ref{prop:2.8} is needed for \cite[Theorem
3.12]{FriRus93}, and our Corollary~\ref{cor:2.2} is needed three
times in \cite[Proposition 4.11]{FriRus93}. Cartan factors are
defined in the next subsection.

\begin{theorem}[\cite{FriRus93},Theorem 8.3] \label{theorem:8.2}
Let $Z$ be an atomic neutral strongly \fss\ satisfying FE, STP,
ERP, and JP. If $Z$ is L-embedded, then
$Z=\oplus_{\alpha}^{\ell^1} J_\alpha$ where each $J_\alpha$ is
isometric to the predual of a Cartan factor of one of the types
1-6. Thus $Z^*$ is isometric to an atomic $JBW^*$-triple. If $Z$
is irreducible, then $Z^*$ is isometric to a Cartan factor.
\end{theorem}

One of our main objectives in this paper is to be able to drop the
assumption of atomicity in this result, i.e. to find a non-ordered
analog of the main theorem of Alfsen-Shultz \cite{AS}. This will
be achieved in our Theorem~\ref{face} below, but at the expense of
some other axioms.

\subsection{\jbst s and ternary rings of
operators}\label{sect:jbst}

A {\it Jordan triple system} is a complex vector space $V$ with a
{\em  triple product} $\{\cdot,\cdot,\cdot\} : V \times V \times V
\longrightarrow V$ which is symmetric and linear in the outer
variables, conjugate linear in the middle variable and satisfies
the Jordan triple identity
\[\{a,b,\{x,y,z\}\} = \{\{a,b,x\},y,z\} - \{x,\{b,a,y\},z\} +
\{x,y,\{a,b,z\}\}. \] A complex Banach space $A$ is called a
$JB^*\mbox{\it -triple}$ if it is a Jordan triple system such that
for each $z\in A,$ the linear map $$ D(z): v\in A\mapsto
\{z,z,v\}\in A $$ is Hermitian, that is, $\|e^{it D(z)}\| = 1$ for
all $t \in \RR$, with non-negative spectrum in the Banach algebra
of operators generated by $D(z)$ and $\Vert D(z)\Vert =\Vert
z\Vert ^2.$ A summary of the basic facts about JB*-triples can be
found in \cite{Russo94} and some of the references therein, such
as \cite{Kaup83},\cite{FriRus85bis}, and \cite{FriRus86bis}.

A $JB^*\mbox{-triple}\; A$ is called a $JBW^*\mbox{\it -triple}$
if it is a dual Banach space, in which case its predual is unique,
denoted by $A_*,$ and the triple product is separately weak*
continuous.  The second dual $A^{**}$ of a $JB^*\mbox{-triple}$ is
a $JBW^*\mbox{-triple.}$

 The $JB^*\mbox{-triples}$ form a large class of Banach spaces which
  include $C^*\mbox{-algebras,}$ Hilbert spaces, spaces of
rectangular matrices, and JB*-algebras. The triple product in a
C*-algebra $\cal A$ is given by
$$ \{x,y,z\} = \,\frac 12\; (xy^*z+ zy^*x). $$
In a JB*-algebra with product $x\circ y$, the triple product is
given by $\{x,y,z\}=(x\circ y^*)\circ z+z\circ (y^*\circ
x)-(x\circ z)\circ y^*$.  An element $e$ in a JB*-triple $A$ is
called a {\it tripotent} if $\{e, e, e \}=e$ in which case the map
$D(e) :  A \longrightarrow A$ has eigenvalues $0,\, {1\over 2}$
and $ 1$, and we have the following decomposition in terms of
eigenspaces $$ A=A_2(e)\oplus A_1(e)\oplus A_0(e) $$ which is
called the {\it Peirce decomposition} of $A$. The ${k\over
2}$-eigenspace $A_k(e)$ is called the {\it Peirce k-space}. The
{\sl Peirce projections} from $A$ onto the Peirce k-spaces are
given by $$ P_2(e) = Q^2(e), \quad P_1(e)= 2(D(e)- Q^2(e)), \quad
P_0(e)= I -2D(e)+Q^2(e) $$ where $Q(e)z = \{ e, z, e\}$ for $z\in
A$. The Peirce projections are contractive.

For any tripotent $v$, the space $A_{2}(v)$ is a JB*-algebra under
the product $x \cdot y =\{x \,\ v \,\ y \}$ and involution
$x^{\sharp}=\{v \,\ x \,\ v \}$.  JBW*-triples have an abundance
of tripotents. In fact, given a JBW*-triple $A$ and $f$ in the
predual $A_*$, there is a unique tripotent $v_f \in A$, called the
{\it support tripotent} of $f$, such that $f \circ P_2(v_f) = f$
and  the restriction $f|_{A_2(v_f)}$ is a {\it faithful positive}
normal functional.

An important class of JBW*-triples are the following six types of
{\it Cartan factors} (see \cite[pp. 292-3]{DF})
:

\begin{description}
\item[type 1] $B(H,K)$,  with triple product $\{x,y,z\}
=\frac{1}{2}(xy^*z +zy^*x),$ \item[type 2] $\{z\in B(H,H): z^t =
-z\},$
\item[type 3]$\{z\in B(H,H): z^t = z\},$ \item[type 4]spin factor (defined below),
\item[type 5]$M_{1,2}({\cal O})$ with triple product
$\{x,y,z\} = \frac{1}{2}(x(y^*z) + z(y^*x)),$\item[type 6]
$M_3({\cal O})$
\end{description}

\noindent where $\cal O$ denotes the 8 dimensional complex
Octonians, $B(H,K)$ is the Banach space of bounded linear
operators between complex Hilbert spaces $H$ and $K$, and $z^t$ is
the transpose of $z$ induced by a conjugation on $H$. Cartan
factors of type 2 and 3 are obviously subtriples of $B(H,H)$, the
latter notation is shortened to $B(H)$, while type 4 can be
embedded as a subtriple of some $B(H)$. The type 3 and 4 are
Jordan algebras with the usual Jordan product $x \circ y =
\frac{1}{2}(xy + yx)$. Abstractly, a {\em spin factor\/} is a
Banach space that is equipped with a complete inner product
$\langle \cdot, \cdot \rangle$ and a conjugation $j$ on the
resulting Hilbert space, with triple product $$ \{ x, y, z\} =
\frac{1}{2} (\langle x,y\rangle z + \langle z,y \rangle x -
\langle  x, jz \rangle jy)$$ such that the given norm and the
Hilbert space norm are equivalent.

An important example of a \jbst\ is a
 {\it ternary ring of operators} (TRO).
This is a subspace of $B(H)$ which is closed under the product
$xy^{\ast}z$. Every TRO is (completely)
 isometric to a corner
$pA(1-p)$ of a C*-algebra $A$. TRO's play an important role in the
theory of quantized Banach spaces (operator spaces), see
\cite{ERuan} for the general theory  and \cite{EOR} for the role
of TRO's. For one thing, as shown by Ruan \cite{Ruan89}, the
injectives in the category of operator spaces are TRO's (corners
of injective C*-algebras) and not, in general, operator algebras.
If $A$ is a TRO and $v$ is a partial isometry in $A$, then
$A_2(v)$ is a \csa\ under the product $(x,y)\mapsto xv^*y$ and
involution $x\mapsto vx^*v$.

 Motivated by a
characterization for JB*-triples  as complex Banach spaces whose
open unit ball is a bounded symmetric domain, we gave in
\cite{NeaRuspj} a holomorphic operator space characterization of
TRO's up to complete isometry. As a consequence, we obtained  a
holomorphic operator space characterization of C*-algebras as
well. Since a closed left ideal in a C*-algebra is a TRO,
Theorem~\ref{theorem:5.5} below will allow us, in our facial
operator space characterization of left ideals (Theorem~\ref{C})
to restrict to TROs from the beginning. The following is the main
result of \cite{NeaRuspj}.

\begin{theorem}[\cite{NeaRuspj},Theorem 4.3]\label{theorem:5.5}
Let $A\subset B(H)$ be an operator space and suppose that
$M_n(A)_0$ is a bounded symmetric domain for some $n \geq 2$. Then
$A$ is n-isometric to a ternary ring of operators (TRO). If
$M_n(A)_0$ is a bounded symmetric domain for all $n\ge 2$, then
$A$ is ternary isomorphic and completely isometric to a TRO.
\end{theorem}

\section{Atomic decomposition of facially symmetric
spaces}\label{sec:atomic}

\subsection{Contractive projections on facially symmetric
spaces}\label{subsect:2.1} In this subsection, we shall assume
that $Z$ is a strongly \fss\ with dual $U=Z^*$. If $\{v_i\}$ is a
countable family of mutually orthogonal \mgt s, then $ v=\sup v_i$
exists as it is the support \gtt\ of $\sum 2^{-i}f_i$, where
$v_i=v_{f_i}$. This fact will be used in the proof of the
following lemma.

\begin{lemma}\label{lem:5.1}
Let $\{v_i\}$ be a countable family of mutually orthogonal \mgt s,
with $ v=:\sup v_i$. Then $v=\sum_i v_i$ (w*-limit).
\end{lemma}
\begin{proof}
 Note first that by \cite[Cor. 3.4(a) and Lemma
1.8]{FriRus89}, for each $n\ge 1$,
\[
\Pi_1^n(P_2(v_i)+P_0(v_i))=\sum_1^n P_2(v_i)+P_0(\sum_1^n v_i),
\]
so by \cite[Cor. 3.4(b)]{FriRus89},
\[
P_2(\sum_1^n v_i)\Pi_1^n(P_2(v_i)+P_0(v_i))=\sum_1^n P_2(v_i),
\]
and hence $\sum_1^n P_2(v_i)$ is a contractive projection. For
$\varphi\in Z$, by orthogonality,
\[
\sum_1^n\|P_2(v_i)\varphi\|=\|\sum_1^nP_2(v_i)\varphi\|\le\|\varphi\|,
\]
so that $\sum_1^\infty\|P_2(v_j)\varphi\|\le\|\varphi\|$ and with
$Q_n:=\sum_1^nP_2(v_i)$ and for $m\ge n$,
\[
\|Q_m\varphi-Q_n\varphi\|=\|Q_n\varphi\|-\|Q_m\varphi\|
\]
so that $Q_n\varphi$ converges to  a limit, call it $Q\varphi$,
and $Q$ is a contractive projection.

For each $x\in U$, $Q_n^*x$ converges in the weak*-topology to
$Q^*x$. Applying this with $x=v$ and recalling that
$P_2(v_i)^*v=v_i$, we obtain $\sum_1^nv_i=Q_n^*v\rightarrow y$ in
the weak*-topology for some $y\in U_2(v)$.  On the other hand,
since $\pair{y}{\sum 2^{-i}f_i}=1$, by \cite[Theorem
4.3(c)]{FriRus89}, we have $F_v\subset F_y$ and therefore by
strong facial symmetry, $y=v+y_0$, where $y_0\in U_0(v)$. Since
$y\in U_2(v)$ we must have $y_0=0$ and hence $y=v$. \qed
\end{proof}

\begin{lemma}\label{lem:5.2} Suppose that $Z$  is neutral and satisfies JP.
Let $\{v_i\}$ be a countable family of mutually orthogonal \mgt s,
with $ v=:\sup v_i$. Then $\cup_{i=1}^\infty [Z_2(v_i)\cup
Z_1(v_i)]$ is norm total in $Z_2(v)+Z_1(v)$.
\end{lemma}
\begin{proof} Let $W$ be the norm closure of the complex span of
$\cup_{i=1}^\infty [Z_2(v_i)\cup Z_1(v_i)]$. We first show that
\begin{equation}\label{eq:W}
Z_2(v)+Z_1(v)\subset W.
\end{equation}
If $\varphi\in Z_2(v)+Z_1(v)$ and $\varphi\not\in W$, then there
exists $x\in U$, $\|x\|\le 1$ with $\pair{x}{\varphi}\ne 0$ and
$\pair{x}{W}=0$. We'll show that $x\in U_0(v)$. Since $\varphi\in
Z_2(v) +Z_1(v)$ implies $\pair{x}{\varphi}=0$, this is a
contradiction, proving (\ref{eq:W}).

Let $s_n=\sum_1^nv_j$ and for $\rho\in Z$ of norm one, let
$\rho=\rho_2+\rho_1+\rho_0$ be its geometric Peirce decomposition
with respect to $s_n$. By JP (for \mgt s), $\rho_2,\rho_1\in W$.
Therefore
$$\pair{s_n\pm
x}{\rho}=\pair{s_n}{\rho_2}\pm\pair{x}{\rho_0}=\pair{P_2(s_n)s_n\pm
P_0(s_n)x}{\rho}$$ so that $$|\pair{s_n\pm
x}{\rho}|\le\|P_2(s_n)s_n\pm
P_0(s_n)x\|=\max(\|P_2(s_n)s_n\|,\|P_0(s_n)x\|)=1.
$$
Thus $\|\sum_1^nv_i\pm x\|=1$ so by Lemma~\ref{lem:5.1}, $\|v\pm
x\|\le 1$. By \cite[Theorem 4.6]{FriRus89}, $v+U_0(v)_1$ is a face
in the unit ball of $U$, and since $v=(v+x)/2+(v-x)/2$, $v\pm x\in
v+U_0(v)_1$, proving $x\in U_0(v)$ and hence (\ref{eq:W}).

To show that equality holds in (\ref{eq:W}), note first that it is
obvious that $Z_2(v_i)\subset Z_2(v)$, and if $\varphi\in
Z_1(v_i)$, then by compatibility,
$P_0(v)\varphi=P_0(v)P_1(v_i)\varphi\in Z_1(v_i)$. But
$P_0(v)\varphi= P_0(v)P_0(v_i)\varphi\in Z_0(v_i)$ so that
$P_0(v)\varphi=0$, as required.\qed\end{proof}

\begin{corollary}\label{cor:5.3}
$P_0(v)=\Pi_{i=1}^\infty P_0(v_i)$ (strong limit).
\end{corollary}
\begin{proof} Let $Q_n=\Pi_1^nP_0(v_i)$. Let $\varphi\in Z$ have geometric
Peirce decomposition $\varphi_2+\varphi_1+\varphi_0$ with respect
to $v$. Since $Z_0(v)\subset Q_n(Z)$,
$Q_n\varphi_0=\varphi_0\rightarrow\varphi_0=P_0(v)\varphi$. It
remains to show that $Q_n(\varphi_2+\varphi_1)\rightarrow 0$. By
Lemma~\ref{lem:5.2}, it suffices to prove that
$Q_n\psi\rightarrow0$ for every $i$ and every $\psi\in
Z_2(v_i)\cup Z_1(v_i)$. But for any $\psi\in Z_k(v_i)$ for
$k=2,1$, $Q_n\psi=Q_nP_k(v_i)\psi=0$ as soon as  $n\ge
i$.\qed\end{proof}

\begin{proposition}
\label{prop:2.8} Suppose that $Z$  is neutral and satisfies JP.
Let $\{u_i\}_{i\in I}$ be an arbitrary family of mutually
orthogonal \mgt s. Then $Q:=\Pi_{i\in I}P_0(u_i)$ exists as a
strong limit and $Q$ is a contractive projection with range
$\cap_{i\in I}Z_0(u_i)$.
\end{proposition}
\begin{proof} Fix $f\in Z$. For each countable set $\lambda\subset I$, let
$g_\lambda:= \Pi_{i\in\lambda}P_0(u_i)f$, which exists as a norm
limit by Corollary~\ref{cor:5.3}.

With $\alpha:=\inf\|g_\lambda\|$, where $\lambda$ runs over the
countable subsets of $I$, we can find a sequence $\lambda_n$ of
countable sets such that $\alpha=\lim \|g_{\lambda_n}\|$, and
hence a
 countable set $\mu=\cup_n\lambda_n\subset I$ such that
$\|g_\mu\|=\alpha$. It remains to prove that
\[
\Pi_{i\in I}P_0(u_i)f=\Pi_{i\in\mu}P_0(u_i)f.
\]
For $\epsilon>0$, choose a finite set $A_0\subset\mu$ such that
for all  finite sets $A$ with $A_0\subset A\subset \mu$,
\[
\|\Pi_{i\in A}P_0(u_i)f-\Pi_{i\in\mu}P_0(u_i)f\|<\epsilon.
\]
By the neutrality of $P_0(u_j)$ and the definition of $\alpha$,
for any $j\not\in\mu$,
\[
P_0(u_j)\Pi_{i\in\mu}P_0(u_i)f=\Pi_{i\in\mu}P_0(u_i)f.
\]
Hence, for any finite subset $B$ with $A_0\subset B\subset I$,
\begin{eqnarray*}\lefteqn{
 \|\Pi_{i\in B}P_0(u_i)f-\Pi_{i\in\mu}P_0(u_i)f\|
}\\
&=&\|\Pi_{i\in B-\mu}P_0(u_i)[\Pi_{i\in
B\cap\mu}P_0(u_i)f-\Pi_{i\in\mu} P_0(u_i)f ]\|\\
&\le&\|\Pi_{i\in B\cap\mu}P_0(u_i)f-\Pi_{i\in\mu} P_0(u_i)f\|
<\epsilon. \qed
\end{eqnarray*}
\end{proof}

\subsection{Jordan decomposition}

In this subsection we introduce the Jordan decomposition property.
We use it in place of atomicity to obtain
Proposition~\ref{prop:38}, which contains the analogs of
\cite[Prop.\ 2.9]{FriRus92} and \cite[Prop.\ 2.4]{FriRus93}.
Lemmas~\ref{prop:1.11}--\ref{lem:yfbr} are taken from an
unpublished note of Yaakov Friedman and the second named author in
1990.
\begin{lemma}\label{prop:1.11}
Let $F$ be a  norm exposed face of the unit ball of a normed space
$Z$, and let $I$ denote the closed unit interval. The following
are equivalent.

(a) $(\mbox{sp}_{{\bf R}}F)_1\subset\mbox{co}(IF\cup-IF)$.

(b) For each non-zero $f\in sp_{\bf R}F, \exists\ g,h\in {\bf
R}^{+}F$ with $f=g-h$ and $\|f\|=\|g\|+\|h\|$.

(c) $\partial (sp_{\bf R}F)_1\subset co(F\cup -F)$.

(d) For each non-zero $f\in sp_{\bf R}F, \exists\ g,h\in {\bf
R}^{+}F$ with $f=g-h$ and $g\perp h$.
\end{lemma}
\begin{proof} (a)$\Rightarrow$(c). If $f\in \mbox{sp}_{{\bf R}}F$ and
$\|f\|=1$, then $f=\lambda\alpha\sigma-(1-\lambda)\beta\tau$, with
$\alpha,\beta, \lambda\in I$ and $\sigma,\tau\in F$. If $\lambda
=0$ or 1, then $f\in\pm F$ so  assume that $0<\lambda<1$. We have
\[
1=\|f\|  =  \|\lambda\alpha\sigma-(1-\lambda)\beta\tau\|
         \leq  \lambda\alpha+(1-\lambda)\beta
        \le \alpha\vee\beta\leq 1.
\]
Since $\lambda<1,\ \alpha=\beta=1.$

(c)$\Rightarrow$(b). If $0\neq f\in sp_{{\bf R}}F$, then
$\|f\|^{-1}f =\lambda\sigma-(1-\lambda)\tau$ with $\lambda\in I$
and $\sigma,\tau\in F$. Since
$\|\lambda\sigma\|+\|(1-\lambda)\tau\|=\lambda+(1-\lambda)=1$, we
have
\[
\|f\|  =  \|f\|(\|\lambda\sigma\|+\|(1-\lambda)\tau\|)
      = \|(\|f\|\lambda)\sigma\|+\|(\|f\|(1-\lambda))\tau\|.
\]

(b)$\Rightarrow$(a). Let $f\in (sp_{{\bf R}}F)_1$ and assume
$0<\|f\|\leq 1$. With $f=g-h$ and $ \|f\|=\|g\|+\|h\|$ with
$g,h\in {\bf R}^{+}F$, we have
\[
f=\|g\|(\|g\|^{-1}g)+\|h\|(-\|h\|^{-1}h) +(1-\|f\|)\cdot 0\in
co(IF\cup -IF).
\]

(d)$\Rightarrow$(b). If $g\perp h$, then
$\|f\|=\|g-h\|=\|g\|+\|h\|$.

(b)$\Rightarrow$(d). If $g,h\in\RR^+F$ and $F=F_x$ for some $x\in
U$ of norm one, then $\|g+h\|=g(x)+h(x)=\|g\|+\|h\|$. Therefore,
$\|g\pm h\|=\|g\|+\|h\|$, {\em i.e.}, $g\perp h$.\qed\end{proof}

\begin{definition}{\rm
A norm exposed face  of the unit ball of a normed space $Z$
satisfies the {\em Jordan decomposition property}  if (one of) the
conditions of Lemma~\ref{prop:1.11} holds.}
\end{definition}

It is elementary that if $F$ satisfies the Jordan decomposition
property, then $\mbox{ext}\, (\mbox{sp}_{{\bf R}}F)_1=\mbox{ext}\,
F\cup \mbox{ext}\, (-F)$.

\begin{lemma}\label{lem:1207}
Let $Z$ be a neutral SFS space, let $F$ be a norm exposed face of
$Z_1$and let $f\in F$. Then $S_{v(f)}(F)\subset F$ and
$P_2(f)F\subset\RR^+F$, where $P_k(f)$ denotes $P_k(v(f))$ for
$k=0,1,2$.
\end{lemma}
\begin{proof} Let $F=F_u$ for some $u\in\gt$. Then by the minimality
property of the polar decomposition (\cite[Theorem
4.3(c)]{FriRus89}), $F_{v(f)}\subset F_u$ and by strong facial
symmetry,  $S_{v(f)}^*u=u$. Thus if $\rho\in F_u,\ \langle
S_{v(f)}\rho,u\rangle=\langle\rho,u\rangle=1$ {\em i.e.},
$S_{v(f)}\rho\in F_u$, which proves the first statement.

By what was just proved,
\[
[P_2(f)+P_0(f)](F)=\frac{I+S_{v(f)}}{2}(F)\subset F.
\]
Thus, if $g\in F$,
\[
\|P_2(f)g\|\left(\frac{P_2(f)g}{\|P_2(f)g\|}\right)+
\|P_0(f)g\|\left(\frac{P_0(f)g}{\|P_0(f)g\|}\right)\in F.
\]
Since $F$ is a face, $P_2(f)g/\|P_2(f)g\|\in F$.\qed\end{proof}

\begin{lemma}\label{prop:1.10}
In a neutral SFS space, the Jordan decomposition is unique
whenever it exists, {\em i.e.}, if $u\in \gt$ and for $i=1,2$, if
$f=\sigma_1-\tau_1=\sigma_2-\tau_2$ with
 $\tau_i,\sigma_i\in\RR^+F_u$,
 $1=\|f\|=\|\sigma_i\|+\|\tau_i\|$,
then $\sigma_1=\sigma_2$ and $\tau_1=\tau_2$.
\end{lemma}
\begin{proof} Apply $P_2(\sigma_1)$ and $P_2(\tau_1)$ to
$f=\sigma_1-\tau_1=\sigma_2-\tau_2$ to obtain $
\sigma_1=P_2(\sigma_1)\sigma_2-P_2(\sigma_1)\tau_2 $ and $
-\tau_1=P_2(\tau_1)\sigma_2-P_2(\tau_1)\tau_2. $ Since
$\sigma_i\perp\tau_i$,
\begin{eqnarray*}
1=\|\sigma_1\|+\|\tau_1\|&\leq&\|P_2(\sigma_1)\sigma_2\|+\|P_2(\sigma_1)
\tau_2\|+\|P_2(\tau_1)\sigma_2\|+\|P_2(\tau_1)\tau_2\|\\
&=&\|[P_2(\sigma_1)+P_2(\tau_1)]\sigma_2\|+\|[P_2(\sigma_1)+P_2(\tau_1)]
\tau_2\|\\
&\leq&\|\sigma_2\|+\|\tau_2\|=1.
\end{eqnarray*}
Therefore
$\|\sigma_1\|=\|P_2(\sigma_1)\sigma_2\|+\|P_2(\sigma_1)\tau_2\|$.

\medskip

Case 1. $P_2(\sigma_1)\tau_2\ne 0$.  In this case,
$P_2(\sigma_1)\sigma_2\ne 0$, otherwise we would have $\sigma_1=0$
and hence $\sigma_1=\sigma_2=0$.  We  then have
\[
\frac{\sigma_1}{\|\sigma_1\|}
=\frac{\|P_2(\sigma_1)\sigma_2\|}{\|\sigma_1\|}\frac{P_2(\sigma_1)\sigma_2}
{\|P_2(\sigma_1)\sigma_2\|}+
\frac{\|P_2(\sigma_1)\tau_2\|}{\|\sigma_1\|}\left(\frac{-P_2(\sigma_1)\tau_2}
{\|P_2(\sigma_1)\tau_2\|}\right).
\]
Since $F_u$ is a face,
\[
\frac{-P_2(\sigma_1)\tau_2}{\|P_2(\sigma_1)\tau_2\|}\in F_u.
\]
On the other hand, by Lemma~\ref{lem:1207},
$P_2(\sigma_1)\tau_2\in\RR^+F_u$, so that $P_2(\sigma_1)\tau_2=0$,
a contradiction, so this case does not occur.

Next, as above, apply $P_2(\sigma_2)$ and $P_2(\tau_2)$ to
$f=\sigma_1-\tau_1=\sigma_2-\tau_2$ to obtain $
\sigma_2=P_2(\sigma_2)\sigma_1-P_2(\sigma_2)\tau_1 $ and $
-\tau_2=P_2(\tau_2)\sigma_1-P_2(\tau_2)\tau_1. $ Since
$\sigma_i\perp\tau_i$, as above we obtain
$\|\sigma_2\|=\|P_2(\sigma_2)\sigma_1\|+\|P_2(\sigma_2)\tau_1\|$.

\medskip

Case 2. $P_2(\sigma_2)\tau_1\ne 0$.  Exactly as in case 1, this
implies that $P_2(\sigma_2)\sigma_1\ne0$ and leads to a
contradiction unless $\sigma_2=0$. So this case does not occur.

\medskip

 Case 3. $P_2(\sigma_2)\tau_1=P_2(\sigma_1)\tau_2= 0$.  In
this case, $\sigma_1=P_2(\sigma_1)\sigma_2$ and
$\sigma_2=P_2(\sigma_2)\sigma_1$, so that
$\|\sigma_1\|=\|\sigma_2\|$. It follows that
$\tau_2(v(\sigma_1))=\pair{P_2(\sigma_1)\tau_2}{v(\sigma_1)}=0$
and
\begin{eqnarray*}
\|\sigma_1\|&=&\sigma_1(v(\sigma_1))=f(v(\sigma_1))=\sigma_2(v(\sigma_1))-\tau_2(v(\sigma_1))\\
&=& \sigma_2(v(\sigma_1))\leq \|\sigma_2\|=\|\sigma_1\|,
\end{eqnarray*}
implying $v(\sigma_2)\leq v(\sigma_1)$. Similarly, using
$P_2(\sigma_2)\tau_1=0$ leads to
$\|\sigma_2\|=\sigma_1(v(\sigma_2))$ and $v(\sigma_1)\leq
v(\sigma_2)$.

Thus $v(\sigma_2)=v(\sigma_1)$, and we now have
\[
\sigma_1=P_2(\sigma_1)f=P_2(\sigma_1)\sigma_2-P_2(\sigma_1)\tau_2=
P_2(\sigma_1)\sigma_2=P_2(\sigma_2)\sigma_2=\sigma_2. \qed
\]
\end{proof}

\begin{lemma}\label{lem:yfbr}
Let $F$ be a norm exposed face satisfying the Jordan decomposition
property. Then
\begin{description}
\item[(a)]$\mbox{sp}_{\RR}F\cap i\mbox{sp}_{\RR}F=\{0\}$.
\item[(b)] If $Z$ is
a neutral strongly symmetric space, then the projection of
$\mbox{sp}_\CC F=\mbox{sp}_\RR F+i\mbox{sp}_\RR F$ onto
$\mbox{sp}_\RR F$ is contractive.
\end{description}
\end{lemma}
\begin{proof} Let $h\in\mbox{sp}_{\RR}F\cap i\mbox{sp}_{\RR}F$, and suppose
that $\|h\|=1$. By Lemma~\ref{prop:1.11},
\[
h=\alpha if+\beta(-ig)=\gamma f_1+\delta (-g_1)
\]
for some $f,g,f_1,g_1\in F$ and $\alpha,\beta,\gamma,\delta\in\RR$
with
 $\alpha\geq
0,\beta=1-\alpha,\gamma\geq 0, \delta=1-\gamma,f\perp g$, and
$f_1\perp g_1$. With $F=F_u$ for some $u\in\gt$, we have
$i(\alpha-\beta)=h(u)=\gamma-\delta$, so that
\begin{equation}\label{eq:*}
h=\frac{1}{2}i(f-g)=\frac{1}{2}(f_1-g_1).
\end{equation}
Applying successively $P_2(f)$ and $P_2(g)$ to (\ref{eq:*}) we
obtain
\[
if=P_2(f)f_1-P_2(f)g_1\mbox{ and }-ig=P_2(g)f_1-P_2(g)g_1.
\]
Therefore
\begin{eqnarray*}
2=\|if\|+\|ig\|&\leq&\|P_2(f)f_1\|+\|P_2(f)g_1\|+\|P_2(g)f_1\|+\|P_2(g)g_1\|\\
&=&(\|P_2(f)f_1\|+\|P_2(g)f_1\|)+(\|P_2(f)g_1\|+\|P_2(g)g_1\|)\\
&\leq&\|f_1\|+\|g_1\|=2.
\end{eqnarray*}
If $P_2(f)f_1=0$, then $if=-P_2(f)f_1\in\RR^+F$, a contradiction.
Similarly, $P_2(f)g_1\ne 0$. Since $iF$ is a face, and
\[
if=\|P_2(f)f_1\|\left(\frac{P_2(f)f_1}{\|P_2(f)f_1\|}\right)
+\|P_2(f)g_1\|\left(\frac{-P_2(f)g_1}{\|P_2(f)g_1\|}\right),
\]
\[
\frac{P_2(f)f_1}{\|P_2(f)f_1\|}\in iF.
\]
On the other hand, by Lemma~\ref{lem:1207},
\[
\frac{P_2(f)f_1}{\|P_2(f)f_1\|}\in F
\]
also. This is a contradiction which proves (a).

Now let $g+ih\in\mbox{sp}_\RR F+i\mbox{sp}_\RR F$. Write $
g=a\rho-b\sigma$ with $\rho\perp\sigma$, $\rho,\sigma\in F$, and
$\|g\|=a+b$. Then $\pair{g}{v_\rho-v_\sigma}=a+b$, and
\begin{eqnarray*}
\|g+ih\|&\ge&|\pair{g+ih}{v_\rho-v_\sigma}|=|a+b+i\pair{h}{v_\rho-v_\sigma}|\\
&=&[(a+b)^2+\pair{h}{v_\rho-v_\sigma}^2]^{1/2}\ge a+b=\|g\|,
\end{eqnarray*}
proving (b). \qed\end{proof}

\medskip

If $Z$ is a dual space, so that each norm exposed face is
weak*-compact, then (b) and the Jordan decomposition property
imply that $\mbox{sp}_\CC F$ is closed, so that
$Z_2(F)=\mbox{sp}_\CC F$.

\begin{definition}
 A WFS space satisfies  property JD if every symmetric
face satisfies the Jordan decomposition property. In this case, we
say that $Z$ is {\em locally base normed}.
\end{definition}
It is important to note that this property is hereditary, that is,
if $Z$ satisfies JD, then so does any geometric Peirce space
$Z_k(u)$.  Indeed, if $F_w\cap Z_k(u)$ is a local face
corresponding to a \gtt\ $w\in U_k(u)$, and $\rho\in
\mbox{sp}_\RR\,[F_w\cap Z_k(u)]$, then $\rho=\alpha g-\beta h$,
with $g,h\in F_w$ and $\|\rho\|=\alpha+\beta$. From this it
follows that $P_k(u)g,P_k(u)h\in F_w$ and $\rho=\alpha
P_k(u)g-\beta P_k(u)h$.

\begin{proposition}\label{prop:38}
Let $Z$ be a locally base normed SFS space.
\begin{description}
\item[(a)] \ii=\mm.
\item[(b)] Suppose furthermore that
$Z$ is  neutral and satisfies JP. Let $v\in \mm$ and suppose that
$w\in \gt$ and $w\top v$. Then $w\in \mm$.
\end{description}
\end{proposition}
 \begin{proof} Let $v\in\ii$ and suppose $F_v$ contains two
distinct elements $f_1,f_2$ and set $f=f_1-f_2$. Then $f=\alpha
g-\beta h$ with $\alpha, \beta\in\RR^+$ and $g,h\in F_v$. By
evaluating at $v$ one sees that $\alpha=\beta=1/2$. Therefore
$F_v$ contains orthogonal elements $g$ and $h$ with orthogonal
supports $v_g$ and $v_h$ such that $v_g\le v$, $v_h\le v$. Since
$v\in \ii$, $v_g=v=v_h$, a contradiction. Thus $F_v$ consists of a
single point and $v\in\mm$. This proves (a).

To prove (b), we first show that $F_w\subset Z_1(v)$. Let
$\psi=\psi_2+\psi_1+\psi_0$ be the Peirce decomposition of
$\psi\in F_w\subset Z_2(w)$ with respect to $v$. We shall show
that $\psi_2=\psi_0=0$.  In the first place, by \cite[Prop.
2.4]{FriRus92}, $\psi_2=P_2(v)\psi=\psi(v)f_v$ and since $v\in
U_1(w)$, $f_v\in Z_1(w)$. (To see this last step, note that for
$k=0,2$,
$P_k(w)f_v=P_k(w)P_2(v)f_v=P_2(v)P_k(w)f_v=\pair{P_k(w)f_v}{v}f_v=
\pair{f_v}{P_k(w)^*v}f_v=0$.) On the other hand, since $v$ and $w$
are compatible, $\psi_2=P_2(v)P_2(w)\psi=P_2(w)P_2(v)\psi\in
Z_2(w)$, showing that $\psi_2\in Z_2(w)\cap Z_1(w)=\{0\}$. Now
$\psi=\psi_1+\psi_0\in F_w$, so $S_v\psi=-\psi_1+\psi_0 \in -F_w$
by \cite[Theorem 2.5]{FriRus89}, so that $\psi_0\in \lc_\RR\,
F_w$. Hence, if
 $\psi_0\ne 0$,we can write $\psi_0/\|\psi_0\|=\lambda\sigma-\mu\tau$ with
$\sigma,\tau\in F_w$, $\lambda,\mu\ge 0$ and $\lambda +\mu=1$.
Since $\sigma,\tau\in F_w$, as shown above, $\sigma_2=\tau_2=0$
and
$\psi_0/\|\psi_0\|=\lambda(\sigma_1+\sigma_0)-\mu(\tau_1+\tau_0)$
so that $\lambda\sigma_1-\mu\tau_1=0$, $\lambda=\mu=1/2$ (since
$\sigma_1(w)=\tau_1(w)=1$) and
$\|\sigma_0-\tau_0\|=\|2\psi_0/\|\psi_0\|\|=2$. Since
\[
2=\|\sigma_0-\tau_0\|\le \|\sigma_0\|+\|\tau_0\|\le 1+1=2,
\]
$\|\sigma_0\|=1$ and by neutrality of $P_0(v)$ (\cite[Lemma
2.1]{FriRus89}, $\sigma_1=0$, implying $\sigma_0=\sigma\in F_w$, a
contradiction as
$\sigma_0(w)=\pair{P_2(w)P_0(v)\sigma}{w}=\pair{P_2(w)\sigma}{P_0(v)^*w}=0$.
Therefore $\psi_0=0$ and $F_w\subset Z_1(v)$.

Now that we know $F_w\subset Z_1(v)$, we show that $F_w$ is a
single point. Suppose to the contrary that there exist $g,h\in
F_w$ with $g\ne h$. Then $f:=g-h$ is a non-zero element of
$\lc_\RR F_w$, so $f=\sigma-\tau$ with $\sigma,\tau\in \RR^+ F_w$
and $\|\sigma\|+\|\tau\|=\|f\|$.  Since $\sigma(w)=\tau(w)$,
$\sigma\ne 0$ and $\tau\ne 0$, and since $\sigma\perp\tau$,
$v_\sigma$ and $v_\tau$ are orthogonal \gtt s in $U_2(w)$ and
hence $U_2(v_\sigma+v_\tau)\subset U_2(w)$. Moreover, by
\cite[Theorems 2.3,2.5]{FriRus89},
\begin{eqnarray*}
U_2(w)&=&\overline{\mbox{sp}}^{w^*}\,\{v_G:G\in\sf,G\subset
Z_2(w)\}\\
&\subset&\overline{\mbox{sp}}^{w^*}\,\{v_G:G\in\sf,G\subset
Z_1(v)\}\\
&\subset&\overline{\mbox{sp}}^{w^*}\,\{v_G:G\in\sf,S_v(G)=-G\}=U_1(v).
\end{eqnarray*}
Then by \cite[Cor. 2.3]{FriRus93}, $v\in U_1(v_\sigma)\cap
U_1(v_\tau)$ and by JP, $v\in U_2(v_\sigma+v_\tau)\subset U_2(w)$,
that is, $w\vdash v$, a contradiction. \qed\end{proof}

\subsection{Rank 2 faces; spin factor}
In this section we assume that $Z$ is a neutral, strongly facially
symmetric, locally base normed space satisfying JP.

\begin{lemma}\label{lem:3.2}
Let $v\in\mm$ and $\varphi\in Z_1(v)$, $\|\varphi\|=1$, and
suppose that $w:=v_\varphi$ is minimal in $U_1(v)$. Then either
$\varphi$ is a global extreme point or the midpoint of two
orthogonal global extreme points.
\end{lemma}
\begin{proof}
 Since $w$ is minimal in $U_1(v)$, the
face $F_w\cap Z_1(v)$ in $Z_1(v)$ exposed by $w$, considered as a
\gtt\ of $U_1(v)$, is the single point $\{\varphi\}$. For any
element $\psi\in F_w$, $P_1(v)\psi=\varphi$, since for $k=0,2$,
$\psi_k(w)=\pair{P_k(v)\psi}{w}=\pair{\psi}{P_k(v)^*w}=0$. Thus
every $\psi\in F_w$ has the form $\psi=\psi_2 +\varphi+ \psi_0$
where $\psi_k=P_k(v)\psi$ for $k=0,2$.

If $F_w=\{\varphi\}$, there is nothing more to prove. So assume
otherwise in the rest of this proof.
 As in the proof of
 Proposition~\ref{prop:38}(a), $F_w$ then contains
two orthogonal elements $\sigma=\sigma_2+\varphi+\sigma_0$ and
$\tau =\tau_2+\varphi+\tau_0$. Further
\begin{eqnarray*}
2&=&\|\sigma-\tau\|=\|\sigma_2-\tau_2\|+\|\sigma_0-\tau_0\|\le
\|\sigma_2\|+\|\tau_2\|+\|\sigma_0\|+\|\tau_0\|\\
&=& \|\sigma_2+\sigma_0\|+\|\tau_2+\tau_0\|\le
\|\sigma\|+\|\tau\|=2.
\end{eqnarray*}

This proves $\|\sigma_2+\sigma_0\|=1=\|\tau_2+\tau_0\|$, and
setting $u:=v_{\sigma_2}\in U_2(v)$ and $\t{u}=v_{\sigma_0}\in
U_0(v)$ one obtains
$\sigma(u+\t{u})=\|\sigma_2\|+\|\sigma_0\|=\|\sigma_2+\sigma_0\|=1$
so $\sigma\in F_w\cap F_{u+\t{u}}$.

We show next that $F_w\cap F_{u+\t{u}}$ is the single point
$\{\sigma\}$.  Suppose to the contrary that $F_w\cap F_{u+\t{u}}$
is not a singleton. Then, as above, it contains two orthogonal
elements $\sigma'$ and $\tau'$ with
$\|\sigma_2'-\tau'_2\|+\|\sigma'_0-\tau_0'\|=2$.

We next claim that
\[
u=v_{\sigma_2'}=v_{\tau'_2}\quad , \quad
\t{u}(\sigma_0')=\|\sigma_0'\| \quad ,\quad
\t{u}(\tau'_0)=\|\tau_0'\|.
\]
Indeed,
\[
1=\pair{u+\t u}{\sigma_2'+\sigma_0'}=\sigma_2'(u)+\sigma_0'(\t u)
\le |\sigma_2'(u)|+|\sigma_0'(\t u)|\le
\|\sigma_2'\|+\|\sigma_0'\|\le 1,
\]
so that $\sigma_0'(\t u)=\|\sigma_0'\|$ and
$\sigma_2'(u)=\|\sigma_2'\|$, and hence $v_{\sigma_2'}\le u$, and
$v_{\sigma_2'}= u$, since $u$, being a multiple of $v$, is
minimal, and hence indecomposable. The proofs for $\tau'_2$ and
$\tau_0'$ are similar.

We next show that there are positive numbers $\lambda,\mu$ and an
extreme point $\rho$ such that $\sigma_2'=\lambda\rho$ and
$\tau_2'=\mu\rho$. Indeed, $\sigma_2'=P_2(v)\sigma'=\sigma'(v)f_v=
\sigma_2'(v)f_v$, where $f_v$ is the extreme point corresponding
to $v\in\mm$, and $\|\sigma_2'\|=\sigma_2'(u)=\sigma_2'(v)f_v(u)$.
Since $u$ is a multiple of $v$, $f_v(u)\ne 0$ and
\[
\sigma_2'=\frac{\|\sigma_2'\|}{f_v(u)} f_v \mbox{ and similarly }
\tau_2'=\frac{\|\tau_2'\|}{f_v(u)} f_v.
\]
Writing $f_v(u)=re^{i\theta}$, we have
$\sigma_2'=\frac{\|\sigma_2'\|}{r} (e^{-i\theta}f_v)$ and
$\tau_2'=\frac{\|\tau_2'\|}{r} (e^{-i\theta}f_v)$.

Finally, assuming without loss of generality that $\lambda\ge
\mu$, we have
\[
2=\|\sigma_2'-\tau_2'\|+\|\sigma_0'-\tau_0'\|=\|\sigma_2'\|-\|\tau_2'\|+\|\sigma_0'-\tau_0'\|\le
1+\|\tau_0'\|-\|\tau_2'\|,
\]
which implies that $\tau_2'=0$ and $\|\tau'_0\|=1$. By neutrality
of $P_0(v)$, $\tau'=\tau_0'$ which is a contradiction.

This proves that  $F_w\cap F_{u+\t{u}}$ is a single point
$\{\sigma\}$ and hence $\sigma=\sigma_2+\varphi+\sigma_0$ is a
global extreme point. Then so is $\t{\sigma}:=-S_v\sigma=
-\sigma_2+\varphi-\sigma_0$ and $\varphi=(\sigma+\t{\sigma})/2$,
completing the proof.\qed\end{proof}

\bigskip

We can now prove  versions of \cite[Prop. 3.2,Lemma 3.6]{FriRus92}
without assuming our space is atomic. First, we need the following
lemma, the conclusion of which is in the hypotheses of \cite[Prop.
3.2,Lemma 3.6]{FriRus92}.

\begin{lemma}\label{lem:3.3}
Let $v\in\mm$, and let $w\in \gt\cap U_1(v)$. Suppose that
$w\not\in\mm$. Then $F_w$ is a rank 2 face, that is, $w=w_1+w_2$
where $w_1$ and $w_2$ are orthogonal \mgt s.
\end{lemma}
\begin{proof} There are two possibilities: (i) $w$ is minimal in $U_1(v)$;
and (ii) $w$ is not minimal in $U_1(v)$.

In case (i), $w$ is the support \gtt\ for some extreme point
$\varphi$ of the unit ball of $Z_1(v)$. Since by assumption
$\varphi$ is not a global extreme point, by Lemma~\ref{lem:3.2},
$\varphi$ is the midpoint of two orthogonal global extreme points,
and therefore $w$ is the sum of two orthogonal \mgt s.

In case (ii), $w=w_1+w_2$ where $w_1,w_2\in\gt\cap U_1(v)$,
$w_1\perp w_2$, and by \cite[Cor. 2.3]{FriRus93}, $w_1\top v$,
$w_2\top v$, so $w_1,w_2\in\mm$ by Proposition~\ref{prop:38}(b).
\qed\end{proof}

\begin{lemma}\label{prop:3.3} Assume that $Z$ also satisfies FE and STP.
Let $v\in\mm$, and let $w\in \gt\cap U_1(v)$. Suppose that
$w\not\in\mm$. Then
\begin{description}
\item[(a)] If $\sigma$ and $\tau$ are orthogonal elements of $F_w$, then
$\sigma$ and $\tau$ are extreme points,
$\sigma+\tau=f_{w_1}+f_{w_2}$ and $v_{\sigma}+v_{\tau}=w$, where
$w=w_1+w_2$ according to Lemma~\ref{lem:3.3}.
\item[(b)] Each norm exposed face of $Z_1$, properly contained in $F_w$,
is a point.
\item[(c)] If $\rho$ is an extreme point of $F_w$, then
there is a unique extreme point $\t \rho$ of $F_w$ orthogonal to
$\rho$.
\item[(d)] With $\xi=(f_{w_1}+f_{w_2})/2$,\ $
F_w=\cup\{[\xi,\rho]: \rho\in \mbox{ext}\, F_w\}$, where
$[\xi,\rho]$ is the line segment connecting $\xi$ and $\rho$.
\end{description}
\end{lemma}
\begin{proof} Case (i). $w$ is minimal in $U_1(v)$.

\medskip

(a) In the proof of Lemma~\ref{lem:3.2}, it was shown that if
$F_w$ contains two orthogonal elements, then these elements are
global extreme points. Once this is known, the equalities
$\sigma+\tau=f_{w_1}+f_{w_2}$ and $v_{\sigma}+v_{\tau}=w$  follow
exactly as in the proof of \cite[Prop. 3.2]{FriRus92}.

(b) Suppose that $F_u\subset F_w$ and $F_u\not=F_w$. By
\cite[Lemma 2.7]{FriRus86} if $\sigma\in F_u$, there exists
$\tau\in F_w$ with $\tau\perp\sigma$. Then $\sigma$ and $\tau$ are
extreme points. Thus $F_u$ consists only of extreme points, and so
it contains only one element.

(c) If $\rho$ is an extreme point of $F_w$, then as in the proof
of (b), there exists an extreme point $\t\rho\in F_w$ orthogonal
to $\rho$. Since by (a), $\rho+\t\rho=f_{w_1}+f_{w_2}$, $\t\rho$
is unique.

(d) The proof is exactly the same as in \cite[Lemma
3.6]{FriRus92}.

\medskip

Case (ii). $w$ is not minimal in $U_1(v)$.

\medskip

In the first place, since $Z_1(v)$ satisfies JD, and $F_w\cap
Z_1(v)$ is not a point, it must contain two orthogonal elements
$g$ and $h$ with orthogonal supports $v_g$ and $v_h$ in $U_1(v)$.
Then by \cite[Cor. 2.3]{FriRus93}, $v_g\top v$, $v_h\top v$, so by
Proposition~\ref{prop:38}(b), $v_g, v_h\in\mm$ and $g,h$ are
global extreme points. After noting that $Z_1(v)$ satisfies FE and
STP (by \cite[Lemma 2.8,Cor.\ 4.12]{FriRus89}), it now follows
exactly as in the proof of case (i) that (a)-(d) hold for the face
$F_w\cap Z_1(v)$. In particular $F_w\cap Z_1(v)=$
 \[
 \{
\lambda\rho+(1-\lambda)\t\rho:\rho,\t\rho \in F_w\cap
Z_1(v)\cap\mbox{ext}\, Z_1,\rho\perp\t\rho, v_\rho+v_{\t\rho}=w,
0\le\lambda\le 1\}.
\]

Now take two orthogonal elements $\sigma,\tau\in F_w$ and Peirce
decompose each one with respect to $v$:
\[
\sigma=\sigma_2+\sigma_1+\sigma_0\quad ,\quad
\tau=\tau_2+\tau_1+\tau_0.
\]
Since $\sigma_1,\tau_1\in F_w$, as noted above we may write
\[
\sigma_1=\lambda\rho+(1-\lambda)\t\rho\quad , \quad
\tau_1=\mu\phi+(1-\mu)\t\phi,
\]
where $\rho$ and $\t\rho$ are orthogonal global extreme points
lying in $F_w\cap Z_1(v)$ with $v_\rho+v_{\t\rho}=w$, and
similarly for $\phi,\t\phi$.

We can partially eliminate $\phi$ and $\t\phi$ as follows. Since
$\tau_1=P_1(v)\tau=P_2(w)P_1(v)\tau\in Z_2(w)$ and
$w=v_\rho+v_{\t\rho}$, by \cite[Lemma 2.3]{FriRus92}
\begin{equation}\label{eq:tau1}
\tau_1=c_1\rho+c_2\t\rho+\psi
\end{equation}
for scalars $c_1,c_2$ and $\psi\in Z_1(\rho)\cap Z_1(\t\rho)$.
Since $|c_1|+|c_2|=
\|c_1\rho+c_2\t\rho\|=\|P_2(v_\rho)\tau_1+P_0(v_{\t\rho})\tau_1\|
\le 1$ and since $1=\tau_1(w)= c_1+c_2+\psi(w)=c_1+c_2$ we have
$c_1+c_2=1$ and $0\le c_1,c_2\le 1$. Denote $c_1$ by $c$ in what
follows.

We shall now prove that
\begin{equation}\label{eq:315}
\tau_0,\sigma_0\in Z_1(v_\rho)\mbox{ and } v_{\tau_0}
v_{\sigma_0}\in\mm,
\end{equation}
and
\begin{equation}\label{eq:316}
\psi\mbox{ in (\ref{eq:tau1}) is zero}.
\end{equation}

To prove (\ref{eq:315}), note that since $v_\rho\in U_1(v)$,
$v_\rho$ is compatible with $v$ , so $P_k(v_\rho)\tau_0\in Z_0(v)$
for $k=2,1,0$. Since $\tau_2=P_2(v)\tau=\pair{\tau}{v}f_v$ and
\begin{eqnarray*}
f_v&=&P_2(\rho)f_v+P_1(\rho)f_v+P_0(\rho)f_v\\
&=&\pair{f_v}{v_\rho}\rho+P_1(\rho)f_v+P_0(\rho)P_2(v)f_v\\
&=&\overline{\pair{\rho}{v}}\rho+P_1(\rho)f_v+\pair{P_0(\rho)f_v}{v}f_v\\
&=&P_1(\rho)f_v,
\end{eqnarray*}
it follows that $\tau_2\in Z_1(\rho)$.  Moreover, since
$S_{v_\rho}^*w=w$, we have $S_{v_\rho}F_w\subset F_w$. Hence
\[
S_{v_\rho}\tau=-\tau_2+c\rho +(1-c)\t\rho-\psi+S_{v_\rho}\tau_0\in
F_w, \]
 and therefore
\[
\frac{\tau+S_{v_\rho}\tau}{2}=c\rho+(1-c)\t\rho+(\tau_0+S_{v_\rho}\tau_0)/2\in
F_w.
\]
Let $\tau':=(\tau_0+S_{v_\rho}\tau_0)/2$. We'll show $\tau'=0$.
Recall that for any $\phi \in Z$,
\[
\|P_{1}(v)\phi+P_{0}(v)\phi\|=\|-S_{v}[P_{1}(v)\phi+P_{0}(v)\phi]\|=
\|P_{1}(v)\phi-P_{0}(v)\phi\|.
\]
Hence, if $\tau'\ne 0$, then
$c\rho+(1-c)\t\rho-(\tau_0+S_{v_\rho}\tau_0)/2\in F_w$, whence
$\tau'\in\ \lc_\RR\, F_w$, and by the property JD,
\[
\tau'/\|\tau'\|=\alpha(\xi_1+\xi_0)-(1-\alpha)(\eta_1+\eta_0)
\]
with $\xi,\eta\in F_w$ and $\alpha\in [0,1]$. Note here that
$$\xi_2=P_2(v)P_2(w)\xi=\pair{P_2(w)\xi}{v}f_v\in Z_1(w),$$ so
$\xi_2=0$ and similarly $\eta_2=0$. As in
Proposition~\ref{prop:38}(b), this implies $\alpha=1/2,
\xi_1=\eta_1$, $\|\xi_0-\eta_0\|=2$ and $\|\xi_0\|=1=\|\eta_0\|$.
By neutrality, $\xi_1=0=\eta_1$, which contradicts the fact that
$\xi=\xi_1+\xi_0\in F_w$.   Thus, $\tau'=0$, proving that
$\tau_0\in Z_1(v_\rho)$. A similar proof shows that $\sigma_0\in
Z_1(v_\rho)$.

Now $v_{\tau_0}\in U_0(v)\cap U_1(v_\rho)$, and if
$v_{\tau_0}\vdash v_\rho$, then $v_\rho\in U_2(v_{\tau_0}) \subset
U_0(v)$, by \cite[Cor.\ 3.4]{FriRus89}, a contradiction. Now by
the two case lemma (\cite[Prop.\ 2.2]{FriRus93}), $v_{\tau_0}\top
v_\rho$ and $v_{\tau_0}$ is a \mgt\ by
Proposition~\ref{prop:38}(b). A similar proof shows that
$\sigma_0\in \mm$. This proves (\ref{eq:315}).

We next prove (\ref{eq:316}). Recall that $
\tau=\tau_2+c\rho+(1-c)\t\rho+\psi+\tau_0, $ and note that
\[
\tau':=-S_vS_{v_\rho}\tau=\tau_2+c\rho+(1-c)\t\rho-\psi+\tau_0,
\]
and $-S_vS_{v_\rho}\sigma=\sigma$. If we let
$\tau'':=(\tau+\tau')/2$ then $\tau,\tau',\tau''\in F_w\cap
\sigma^{\perp}$, in particular $\psi=\tau''-\tau'\in
\sigma^\perp$. Suppose that $\psi$ is not a multiple of a global
extreme point. Since $\psi\in Z_1(v_\rho)\cap Z_1(v_{\t\rho})$,
and $v_\psi$ is not minimal, we have $v_\rho,v_{\t\rho}\in
U_2(v_\psi)$ and $w=v_\rho+v_{\t\rho}\in U_2(v_\psi)$.
 But $v_\sigma\le w\in
U_2(v_\psi)$ and $v_\sigma\perp v_\psi$, implying $v_\sigma\in
U_0(v_\psi)\cap U_2(v_\psi)$, a  contradiction.

We conclude that  $\psi=\alpha\varphi$ is a multiple of a global
extreme point $\varphi$. From (\ref{eq:tau1}), if $\alpha\ne 0$,
then $\varphi$ is a difference of two elements of $F_w$, hence an
extreme point of $(\lc_\RR\, F_w)_1$, which implies that
$\varphi\in F_w\cup F_{-w}$. This is a contradiction since
$\pm\alpha =\alpha\varphi(w)=\psi(w)=\psi(v_\rho+v_{\t\rho})=0$.
Hence $\alpha=0$ proving (\ref{eq:316}).

We next show that
$F_w\cap\{\sigma\}^\perp\cap\{\tau\}^\perp=\emptyset$. Suppose
there exists a point $\tau'$ lying in
$F_w\cap\{\sigma\}^\perp\cap\{\tau\}^\perp$. By the above
calculations, one member of the set
$\{\tau_{1},(\tau')_{1},\sigma_{1}\}$ is a convex combination of
the other two. From this it follows exactly as in the proof of
Lemma~\ref{lem:3.2} that the corresponding convex combination of
two elements of the orthogonal set $\{\tau,\tau',\sigma\}$ is an
extreme point, which is a contradiction. Thus
$F_w\cap\{\sigma\}^\perp\cap\{\tau\}^\perp=\emptyset$.

We can now complete the proof of (a), and (b)-(d) will follow as
in case (i). If $F_{v_\tau}\ne\{\tau\}$, then by JD, $F_{v_\tau}$
contains two orthogonal elements $g,h$. But we have proved that in
this case $F_w\cap\{g\}^\perp\cap\{h\}^\perp=\emptyset$. However,
this set contains $\sigma$ and this contradiction shows that
$\tau$ (and by symmetry $\sigma$) is an extreme point. This
completes the proof of Lemma~\ref{prop:3.3}.\qed\end{proof}

\medskip

Once we know the result of Lemma~\ref{prop:3.3} above, the proof
in \cite{FriRus92} shows that the main result of \cite{FriRus92}
holds with atomic replaced by JD and JP. We formalize this in the
next proposition.

\begin{proposition}\label{theorem:3.6}
Let $Z$ be a neutral strongly facially symmetric locally base normed
space  which satisfies FE, STP,  and JP. If $v\in\mm$ and $u\in \gt\cap
U_1(v)$, then $Z_2(u)$ is isometric to the dual of a complex spin
factor.
\end{proposition}
\begin{proof} The argument in \cite{FriRus92}, from \cite[Corollary
3.7]{FriRus92} to \cite[Theorem 4.16]{FriRus92}
 uses only the following results from \cite{FriRus92} and
does not otherwise invoke the atomic assumption made there:
\cite[Prop.\ 2.9,Cor.\ 2.11,Prop.\ 3.2,Lemma 3.6]{FriRus92}.

On the one hand, \cite[Prop.\ 2.9]{FriRus92} and \cite[Cor.\
2.11]{FriRus92} remain true if atomic is replaced there by JD and
JP, as shown in our Proposition~\ref{prop:38}(a). On the other
hand, \cite[Prop\ 3.2]{FriRus92} remains true if atomic is
replaced by JD and JP, as shown in our
Lemma~\ref{prop:3.3}(a),(b),(c); and \cite[Lemma 3.6]{FriRus92}
remains true if atomic is replaced by JD and JP, as shown in our
Lemma~\ref{prop:3.3}(d). Thus Proposition~\ref{theorem:3.6} is
proved. \qed\end{proof}

\subsection{Atomic decomposition}
The following is the main result of this section.

\begin{theorem}\label{atdec}
Let $Z$ be a locally base normed neutral strongly \fss\ satisfying the
pure state properties and JP.
  Then $Z=Z_a\oplus^{\ell^1}N$,
where $Z_a$ and $N$ are strongly \fss s satisfying the same
properties as $Z$, $N$ has no extreme points in its unit ball, and
$Z_a$ is the norm closed complex span of the extreme points of its
unit ball.
\end{theorem}
\begin{proof} If $Z$ has no extreme points in its unit ball, there is
nothing to prove. If it has an extreme point, then there exists a
maximal family
 $\{u_i\}_{i\in I}$ of mutually orthogonal \mgt s.
Let $Q:=\Pi_{i\in I}P_0(u_i)$ be the contractive projection on $Z$
with $Q(Z)=\cap_{i\in I}Z_0(u_i)$ guaranteed by
Proposition~\ref{prop:2.8}. We shall show that $N:=Q(Z)$ and
$Z_a:=(I-Q)(Z)$ have the required properties. By maximality, $N$
has no extreme points in its unit ball.

For  a finite subset $A$ of $I$ and $Q_A:= \Pi_{i\in A}P_0(u_i)$,
by JP,
\begin{eqnarray}\label{eq:399}\lefteqn{
(I-Q_A)(Z)=Z_2(\Sigma_A u_i)\oplus Z_1(\Sigma_Au_i)}\\\nonumber
 &=&\left(\oplus_A
Z_2(u_i)\right)\oplus\left(\oplus_{i\ne j}[Z_1(u_i)\cap
Z_1(u_j)]\right)\oplus\left(\oplus_A [Z_1(u_i)\cap
Z_0(\Sigma_{j\ne i}u_j)]\right).
\end{eqnarray}
Since $I-Q_A\rightarrow I-Q$ strongly, it follows that every
element of $(I-Q)(Z)$ is the norm limit of elements from
$\cup_A(I-Q_A)(Z)$. Since obviously $Z_2(u_i)\perp Q(Z)$, in order
to prove $Z_a\perp N$, it suffices to prove that for every $i\in
I$,
\begin{equation}\label{eq:z1}
Z_1(u_i)\perp Q(Z).
\end{equation}

For each $i$, let $Q_i=\Pi_{j\in I-\{i\}}P_0(u_j)$ and for
 $\varphi\in Z_1(u_i)$, write $\varphi=Q_i\varphi+(I-Q_i)\varphi$.
Note that
\[
Q_i(Z_1(u_i))=Z_1(u_i)\cap[\cap_{j\in I-\{i\}}Z_0(u_j)]
\]
and that
\[
(I-Q_i)(Z_1(u_i))\mbox{ is the norm closure of
}\oplus^{\mbox{finite}}_{j\in I-\{i\}}[Z_1(u_i)\cap Z_1(u_j)].
\]
For the latter, note that for a finite subset $A\subset I-\{i\}$,
if $Q_{i,A}$ denotes the partial product for $Q_i$, then
\begin{eqnarray*}
(I-Q_{i,A})P_1(u_i) &=&\sum_A P_2(u_j)P_1(u_i)+\sum_{k\ne
l}P_1(u_k)P_1(u_l)
P_1(u_i)\\
&+&\sum_AP_1(u_j)P_0(\sum_{k\ne j}u_k)P_1(u_i)\\
&=&0+0+\sum_AP_1(u_j)P_1(u_i).
\end{eqnarray*}
Thus, $(I-Q_i)\varphi$ can be approximated in the norm by elements
from spaces of the form $\oplus_{j\in A} [Z_1(u_i)\cap Z_1(u_j)]$,
where $A$ is a finite subset of $I-\{i\}$.

Now (\ref{eq:z1}) is reduced to proving that $Q_i(Z_1(u_i))\perp
Q(Z)$ and $(I-Q_i)(Z_1(u_i))\perp Q(Z)$. Since $Z_1(u_i)\cap
Z_1(u_j)\subset Z_2(u_i+u_j)$ and $Q(Z)\subset Z_0(u_i+u_j)$, it
is clear that $[Z_1(u_i)\cap Z_1(u_j)]\perp Q(Z)$. It remains to
show that
\begin{equation}\label{eq:z2}
\left(Z_1(u_i)\cap [\cap_{j\in I-\{i\}}Z_0(u_j)]\right)\perp Q(Z).
\end{equation}

Suppose $g\in Z_1(u_i)\cap [\cap_{j\in I-\{i\}}Z_0(u_j)]$ and
$h\in Q(Z)$. Then either $v_g\vdash u_i$ or $v_g\top u_i$.  In the
first case, since by Proposition~\ref{theorem:3.6}, $U_2(v_g)$ is
isometric to a spin factor, there is a \mgt\ $\t u_i$ with $\t
u_i\perp u_i$ and $\t u_i\in U_0(\Sigma_ {j\in I-\{i\}}u_j)$. This
contradicts the maximality. Therefore $v_g$ is a \mgt\ and $g$ is
a multiple of an extreme point $\psi$. If $h=h_2+h_1+h_0$ is the
geometric Peirce decomposition of $h$ with respect to $v_g$, then
since $v_g$ is compatible with all the $u_k$, $h_j\in Q(Z)$. Now
$h_2$ is also a multiple of $\psi$ and $\psi\in Z_1(u_i)$; hence
$h_2\in Z_0(u_i)\cap Z_1(u_i)=\{0\}$. Since $v_{h_1}\in U_1(v_g)$,
either $v_{h_1}\vdash v_g$ or $v_{h_1}\top v_g$. In the first case
we would have $v_g\in U_2(v_{h_1})\subset Q(Z)$, a contradiction.
In the second case, $h_1$ would be a multiple of $\psi$, again a
contradiction. We conclude that $h_1=0$ and therefore $h=h_0\in
Z_0(v_g)$ so that $g\perp h$ as required, proving (\ref{eq:z2})
and thus the decomposition $Z=Z_a\oplus^{\ell^1}N$.

It is elementary that all the properties of $Z$ transfer to any
$L$-summand. Finally, the set of extreme points of the unit ball
of $Z$ which lie in  $(I-Q)(Z)$ are norm total in $(I-Q)(Z)$,
since every element from the right side of (\ref{eq:399}) is a
linear combination of at most two extreme points by
Lemmas~\ref{lem:3.2} and ~\ref{prop:3.3}(d). \qed\end{proof}
\section{Characterization of one-sided ideals in \csa
s}\label{sec:ideal}
\subsection{Contractive projections on Banach spaces}\label{sec:2}

An interesting question about general Banach spaces, which is
relevant to this paper,
 is to determine under what conditions  the
intersection of 1-complemented subspaces is itself 1-complemented.
Although  this may be true if the contractive projections onto the
subspaces form a commuting family, we have been unable to prove it
or find it in the literature, without adding some other
assumptions. The hypothesis of weak sequential completeness used
in Corollary~\ref{cor:7.2} and  Proposition~\ref{theorem:2.4} is
satisfied in $L$-embedded spaces, as noted in
subsection~\ref{sect:fss}.

\begin{lemma}
Let $X$ be a Banach space and let $\{P_i\}_{i\in I}$ be a family
of commuting contractive projections on $X$. Then $W:=\cap_{i\in
I}P_i^*(X^*)$ is the range of a contractive projection on $X^*$.
\end{lemma}
\begin{proof} Let $\ff$ denote the collection of finite subsets of $I$. For
each $A\in\ff$, let $Q_A=\Pi_{i\in A}P_i$. Since the unit ball
$B(X^*)_1$ is compact in the weak*-operator topology (=
point-weak*-topology), there is a subnet
$\{R_\delta\}_{\delta\in\dd}$ of the net $\{Q_A^*\}_{A\in\ff}$
converging in this topology to an element $R\in B(X^*)_1$. Thus
$R_\delta=Q^*_{u(\delta)}$, where $u:\dd \rightarrow\ff$ is a
finalizing map ($\forall A\in\ff,\exists \delta_0 \in\dd,
u(\delta)\ge A,\forall \delta\ge \delta_0$), and for every $x\in
X^*$ and $f\in X$,
\[
\pair{Rx}{f}=\lim_\delta\pair{R_\delta x}{f}.
\]

It is now elementary to show that $R^2=R$ and $Rx=x$ if and only
if $x\in W$. For completeness, we include the details.

For $x\in X^*,\ f\in X$,
\begin{eqnarray*}
\pair{R^2x}{f}&=&\lim_\delta\pair{R_\delta Rx}{f}=
\lim_\delta\pair{Rx}{Q_{u(\delta)}f}\\
&=&\lim_\delta\lim_{\delta'}\pair{R_{\delta'} x}{Q_{u(\delta)}f}=
\lim_\delta\lim_{\delta'}\pair{x}{Q_{u(\delta')}Q_{u(\delta)}f}\\
&=&\lim_{\delta'}\pair{x}{Q_{u(\delta')}f}=
\lim_{\delta'}\pair{R_{\delta'}x}{f}=\pair{Rx}{f}.
\end{eqnarray*}
Thus $R^2=R$.

If $x\in W$, then $Q_A^*x=x$ for every $A\in\ff$, so that
$\pair{Rx}{f}=\lim_\delta\pair{R_\delta x}{f}=
\lim_\delta\pair{Q_{u(\delta)}^*x}{f}=\pair{x}{f}$, so that
$Rx=x$.

Conversely, if $Rx=x$, then
\begin{eqnarray*}
\pair{P_i^*x}{f}&=&\pair{P_i^*Rx}{f}=\lim_\delta\pair{R_\delta
x}{P_if}\\
&=&\lim_\delta\pair{P_i^*Q_{u(\delta)}^*x}{f}=\lim_\delta\pair{R_\delta
x}{f} =\pair{Rx}{f}=\pair{x}{f},
\end{eqnarray*}
so that $x\in W$. \qed\end{proof}

\bigskip

We cannot conclude from the above proof that $\cap_{i\in I}P_i(X)$
is the range of a contractive projection on $X$. On the other
hand, we have the following two immediate consequences.

\begin{corollary}\label{cor:2.2}
Let $X$ be a reflexive Banach space, and let $\{P_i\}_{i\in I}$ be
a family of commuting contractive projections on $X$ with ranges
$X_i=P_i(X)$. Then $Y:=\cap_{i\in I}X_i$ is the range of a
contractive projection on $X$.
\end{corollary}

\begin{corollary}\label{cor:7.2}
Let $X$ be a weakly sequentially complete Banach space, and let
$\{P_i\}_{i\in \NN}$ be a  sequence of commuting contractive
projections on $X$ with ranges $X_i=P_i(X)$. Then $Y:=\cap_{i\in
\NN}X_i$ is the range of a contractive projection on $X$.
\end{corollary}
\begin{proof} With $Q_n=P_1\cdots P_n$, there is a subsequence $Q_{n_k}^*$
converging to an element $R\in B(X^*)_1$ in the weak*-operator
topology, that is, for $x\in X^*$ and $f\in X$,
$\pair{x}{Q_{n_k}f}\rightarrow\pair{Rx}{f}$, so that
$\{Q_{n_k}f\}$ is a weakly Cauchy sequence. By assumption,
$Q_{n_k}f$ converges weakly to an element $Sf$, and it is
elementary to show that $R=S^*$, and $S$ is a contractive
projection on $X$ with range $Y$. \qed\end{proof}

\begin{proposition}\label{theorem:2.4}
Let $X$ be a weakly sequentially complete Banach space, and let
$\{P_i\}_{i\in I}$ be a family of neutral commuting contractive
projections on $X$ with ranges $X_i=P_i(X)$. Then $Y:=\cap_{i\in
I}X_i$ is the range of a contractive projection on $X$.
\end{proposition}
\begin{proof} We note first that for any countable subset $\lambda\subset
I$, by Corollary~\ref{cor:7.2}, there is a contractive projection
$Q_\lambda$ (not necessarily unique).
 with range $\cap_{i\in\lambda}X_i$
Now, for $f\in X$, define
\[
\alpha_f=\inf_\lambda\inf_{Q_\lambda} \|Q_\lambda f\|.
\]
There exists a sequence $\lambda^{(n)}$ and a choice of
contractive projection $Q_{\lambda^{(n)}}$ with
$\alpha_f\le\|Q_{\lambda^{(n)}}f\|\le\alpha_f+1/n$. Set
$\mu=\cup_n\lambda^{(n)}$ and let $Q_\mu$ be a contractive
projection on $X$ with range $\cap_{i\in\mu}X_i$. Since
$Q_\mu(X)\subset
 Q_{\lambda^{(n)}}(X)$, we have $\|Q_\mu f\|=\|Q_\mu
Q_{\lambda^{(n)}}f\| \le \|Q_{\lambda^{(n)}}f\|$ implying
$\alpha_f=\|Q_\mu f\|$, and so $\|Q_\mu f\|\le \|Q_\lambda f\|$
for all countable subsets $\lambda$ of $I$.

If $Q'_\mu$ is any other contractive projection with range
$Q_\mu(X)$, then $Q_\mu f=Q'_\mu Q_\mu f=Q_\mu Q'_\mu f=Q'_\mu f$
so that we may unambiguously define an element $Qf\in\cap_{i\in
\mu}X_i$ by $Qf:=Q_\mu f$. By the neutrality of the projections,
it follows that $Qf\in\cap_{i\in I}X_i$. Indeed, if $j\in I-\mu$,
then
\[
\|Q_{\mu\cup\{j\}}f\|=\|P_jQ_\mu f\|\le \|Q_\mu f\|\le
\|Q_{\mu\cup\{j\}}f\|,
\]
and by the neutrality of $P_j$, $P_jQ_\mu f=Q_\mu f$. Hence
$Qf\in\cap _{i\in I}X_i$. Conversely, if $f\in \cap_{i\in I}X_i$,
then in particular, $f\in Q_\mu(X)$, so $Qf=Q_\mu f=f$.

We have shown that $Q$ is a nonlinear nonexpansive projection of
$X$ onto $Y$. It remains to show that $Q$ is actually linear. For
this it suffices to observe that, by neutrality, if $Qf=Q_\mu f$,
then $Qf=Q_\lambda f$ for any countable set $\lambda\supset \mu$.
Then, if $f,g\in X$ and $Qf=Q_\mu f,\ Qg=Q_\nu g$, and
$Q(f+g)=Q_\sigma(f+g)$ for suitable countable sets
$\mu,\nu,\sigma$ of $I$, then with $\tau=\mu\cup \nu\cup\sigma$,
\[
Q(f+g)=Q_\tau(f+g)=Q_\tau f+Q_\tau g=Qf+Qg.\qed
\]
\end{proof}

\subsection{Characterization of predual of Cartan factor}
In this subsection we show that the entire machinery of
\cite{FriRus93} can be repeated with appropriate modifications to
yield a variation of the main result of \cite{FriRus93} to
non-atomic \fss s satisfying JD, and stated in
Proposition~\ref{theorem:3.7} below.  As noted below, the
assumption that $Z$ is $L$-embedded in its second dual needs to be
added to the assumptions in \cite{FriRus93}. As was done in the
proof of Propositiion~\ref{theorem:3.6}, we shall explicitly
indicate the modifications needed in \cite{FriRus93}, section by
section, to prove Proposition~\ref{theorem:3.7}.

In the proof of \cite[Lemma 1.2]{FriRus93} it was stated that the
intersection of a certain family of 1-complemented subspaces, is
itself 1-complemented. As noted in Section~\ref{sec:2}, this is
problematical in general. However, \cite[Lemma 1.2]{FriRus93} is
used in \cite{FriRus93} only in the context of a reflexive Banach
space, hence it is covered by Corollary~\ref{cor:2.2}. The role of
the assumption of atomic in \cite[Proposition 1.5]{FriRus93} is to
obtain the property expressed in Proposition~\ref{prop:1.1}(b).
But as shown in Theorem~\ref{atdec}, this property will be
available to us. Thus \cite[Section 1]{FriRus93} is valid with
atomic replaced by JD.

By Proposition~\ref{prop:38}(b) and Lemma~\ref{lem:3.3}
respectively, \cite[Proposition 2.4]{FriRus93} and
\cite[Proposition 2.5]{FriRus93} remain true with atomic replaced
by JD. \cite[Corollary 2.7]{FriRus93} depends only on
\cite[Proposition 2.5]{FriRus93} and
Proposition~\ref{theorem:3.6}, while the part of \cite[Lemma
2.8]{FriRus93} concerned with the property FE
 is immediate from
\cite[Corollary 2.7]{FriRus93} and \cite[Proposition
2.4]{FriRus93}. Finally, \cite[Corollary 2.9]{FriRus93} is
immediate from \cite[Proposition 2.9]{FriRus92} which, as already
remarked in the proof of Proposition~\ref{theorem:3.6}, remains
true with
 atomic replaced by JD (Proposition~\ref{prop:38}(a)).
Thus \cite[Section 2]{FriRus93} is valid with atomic replaced by
JD.

The only reliance on atomicity in \cite[Section 3]{FriRus93}
occurs in \cite[Lemma 3.2]{FriRus93} and \cite[Proposition
3.7]{FriRus93}. The former depends only on \cite[Corollaries 2.7
and 2.9]{FriRus93} and the latter on \cite[Proposition
1.5]{FriRus93}, which as just noted, are both valid with atomic
replaced by JD. In the proof of \cite[Proposition 3.12]{FriRus93}
it was stated that the intersection of a  family of Peirce-0
subspaces of an orthogonal family of \mgt s is 1-complemented, and
in fact the net of partial products converges strongly to the
projection on the intersection. As no proof was provided for this
in \cite{FriRus93}, we provided a proof in
Proposition~\ref{prop:2.8}. Recall that Proposition~\ref{prop:2.8}
was also the key engredient of the proof of the atomic
decomposition in Theorem~\ref{atdec} above.

With these remarks  we can now assert the following modification
of \cite[Theorem 3.14]{FriRus93}, the main result of \cite[Section
3]{FriRus93}.

\begin{lemma}
Let $Z$ be a neutral locally base normed SFS space  and assume the pure
state properties FE,ERP, and STP, and the property JP. Assume
that there exists a minimal geometric tripotent $v$ with $U_1(v)$
of rank 1 and a geometric tripotent $u$ with $u\vdash v$. Then $U$
has an $M$-summand which is linearly isometric with the complex
\jbwst\ of all symmetric ``matrices'' on a complex Hilbert space
(Cartan factor of type 3). In particular, if $Z$ is irreducible,
then $Z^*$ is isometric to a Cartan factor of type 3.
\end{lemma}

The only possible reliance on atomicity in \cite[Section
4]{FriRus93} occurs in \cite[Lemma 4.9]{FriRus93} and
\cite[Proposition 4.11]{FriRus93}. The former depends only on
\cite[Lemma 3.6]{FriRus92}, which is valid in the presence of JD
by Lemma~\ref{prop:3.3}(d), and the latter on \cite[Lemma
1.2]{FriRus93}, which as noted above is needed only for reflexive
Banach spaces. However, \cite[Lemma 4.9]{FriRus93} states
explicitly that reflexivity. Note that the ``classification
scheme'', embodied in \cite[Proposition 4.20]{FriRus93} does not
involve atomic so is valid in the presence of JD.

The only reliance on atomicity in \cite[Section 5]{FriRus93}
occurs in \cite[Lemma 5.2]{FriRus93}, which depends on
\cite[Corollary 2.11]{FriRus92}. As already noted, the latter is
valid in the presence of JD. In the proof of \cite[Lemma
5.5]{FriRus93} it was stated that the intersection of a  family of
Peirce-0 subspaces of a family of geometric tripotents which are
either orthogonal or collinear is 1-complemented, and in fact the
net of partial products converges strongly to the projection on
the intersection. As no proof was provided for this in
\cite{FriRus93}, we provided a proof of the 1-complementedness of
the intersection in Proposition~\ref{theorem:2.4}. This is the
only place in this paper and one of two places in \cite{FriRus93}
where the assumption of L-embeddedness is used.
 Although it is problematical whether
 the strong convergence of the partial products
exists, nevertheless, it is sufficient to take a subnet of the net
of partial sums in the proof of \cite[Lemma 5.5]{FriRus93}. The
same remark applies to \cite[Lemma 6.6]{FriRus93}.

With these remarks  we can now assert the following modification
of \cite[Theorem 5.10]{FriRus93}, the main result of \cite[Section
5]{FriRus93}.

\begin{lemma}
Let $Z$ be a neutral locally base normed SFS space of spin degree 4, which
is L-embedded and which satisfies FE,STP,ERP, and JP. Then $Z$
has an $L$-summand which is linearly isometric to the predual of a
Cartan factor of type 1. In particular, if $Z$ is irreducible,
then $Z^*$ is isometric to a Cartan factor of type 1.
\end{lemma}

The only  reliance on atomicity in \cite[Section 6]{FriRus93}
occurs in \cite[Lemma 6.2]{FriRus93}. However, this dependence is
on earlier results which have been established in the presence of
JD. As noted above for \cite[Lemma 5.5]{FriRus93}, \cite[Lemma
6.6]{FriRus93} holds under the assumption of $L$-embeddedness.

With these remarks  we can now assert the following modification
of \cite[Theorem 6.8]{FriRus93}, the main result of \cite[Section
6]{FriRus93}.

\begin{lemma}
Let $Z$ be a neutral locally base normed SFS space of spin degree 6, which
is L-embedded and which satisfies FE,STP,ERP, and JP. Then $Z$
has an $L$-summand which is linearly isometric to the predual of a
Cartan factor of type 2. In particular, if $Z$ is irreducible,
then $Z^*$ is isometric to a Cartan factor of type 2.
\end{lemma}

The results of \cite[Sections 7 and 8]{FriRus93} carry over
verbatim in the presence of JD.  The proof of \cite[Theorem
7.1]{FriRus93} on pages 75--79 of \cite{FriRus93} yields the
following modification.

\begin{lemma}
Let $Z$ be a neutral locally base normed SFS space which satisfies FE,
STP, ERP, and JP, and let $v,\t v$ be orthogonal minimal geometric
tripotents in $U:=Z^*$ such that the dimension of $U_2(v+\t v)$ is
8 and $U_1(v+\t v)\ne\{0\}$. Then there is an $L$-summand of $Z$
which is isometric to the predual of a  Cartan factor of type 5,
\ie, the 16 dimensional \jbwst\ of 1 by 2 matrices over the
Octonions. In particular, if $Z$ is irreducible, then $Z^*$ is
isometric to the Cartan factor of type 5.
\end{lemma}

Similarly, the proof of \cite[Theorem 7.8]{FriRus93} appearing on
pages 79--82 of \cite{FriRus93} yields the following modification.

\begin{lemma}
Let $Z$ be a neutral locally base normed SFS space of spin degree 10 which
satisfies FE,STP,ERP, and JP, and has no $L$-summand of type $I_2$.
Then $Z$ contains an $L$-summand which is isometric to the predual
of a Cartan factor of type 6, \ie, the 27 dimensional \jbwst\ of
all 3 by 3 hermitian matrices over the Octonions. In particular,
if $Z$ is irreducible, then $Z^*$ is isometric to the Cartan
factor of type 6.
\end{lemma}

Finally, the proof of  \cite[Theorem 8.2]{FriRus93} on pages
83--84 of \cite{FriRus93} yields the following modification.

\begin{proposition}\label{theorem:3.7}
Let $Z$ be a neutral locally base normed strongly \fss\ satisfying FE,
STP, ERP, which is L-embedded and which satisfies JP. For any
\mgt\ $v$ in
$U$, there is an $L$-summand $J(v)$ of $Z$ isometric to the
predual of a Cartan factor of one of the types 1-6 such that $\hat
v\in J(v)$. If $Z$ is the norm closure of the complex linear span
of its extreme points, then it is isometric to the predual of an
atomic \jbwst.
\end{proposition}

\subsection{Spectral duality and Characterization of dual ball of \jbst}

 If $Z$ is an L-embedded, locally base normed, neutral strongly
facially symmetric space satisfying JP and the pure state
properties, then by Proposition~\ref{theorem:3.7} and
Theorem~\ref{atdec}, its dual $Z^*$ is a direct sum
$Z^*=(Z_a)^*\oplus^{\ell^\infty}N^*$ where $(Z_a)^*$ is isometric
to an atomic \jbwst. We shall identify $(Z_a)^*$ with this \jbwst\
in what follows.

\begin{lemma}\label{lem:316}
Suppose that $Z$ is as above and assume that $Z$ is the dual of a
Banach space $B$. For $a\in B$, if $\hat a$ denotes the canonical
image of $a$ in $Z^*$, and $Q$ is the projection of $Z^*$ onto
$(Z_a)^*$, then $\|Q\hat a\|=\|a\|$.
\end{lemma}
\begin{proof} For $a\in B$ with $\|a\|=1$, let $g$ be an extreme point of
the nonempty convex w*-compact set $\{f\in Z:\|f\|=1=f(a)\}$. Then
$g\in\mbox{ ext}\, Z_1$, so $g$ vanishes on $N^*$. Thus
\[
1=\|a\|=\|\hat a\|\ge \|Q\hat a\|\ge|\pair{Q\hat
a}{g}|=|\pair{\hat a}{g}|=|\pair{g}{a}|=1. \qed
\]
\end{proof}

\medskip

In order to show that the space $B$ is isometric to a \jbst, it
suffices to show that the image of the map $a\mapsto Q\hat a$ is
closed under the cubing operation in $(Z_a)^*$, and is hence a
subtriple of  $(Z_a)^*$. To show this we need a spectral
assumption on the elements of $B$. To make this definition, we
need a lemma.

\begin{lemma}\label{lem:spectral}
Let $Z$ be a neutral WFS space satisfying PE.
 Let $\{F_B:B\in {\cal B}\}$ be a family of
norm closed faces of $Z_1$, where $\cal B$ denotes the set of
non-empty Borel subsets of the closed interval $[a,b]$.

\begin{description}
\item[(a)] Suppose that
\begin{description}
\item[(i)] if $B_1\cap B_2=\emptyset$, then
$F_{B_1}\perp F_{B_2}$ and $v_{B_1\cup B_2}=v_{B_1}+v_{B_2}$.
\end{description}

For $f\in C[a,b]$, if $P=\{s_0,\ldots,s_n\}$ is a partition  of
$[a,b]$ and $T=\{t_1,\ldots,t_n\}$ are points  with $s_{i-1}\le
t_i\le s_i$, the Riemann sums $S(P,T,f)=\sum_1^n
f(t_j)v_{(s_{j-1},s_j]}$ converge in norm to an element $\int f\,
dv_B=\int f(t)\, dv_B(t)$ of $Z^*$ as the mesh
$|P|=\min\{s_j-s_{j-1}\}\rightarrow 0$.

\item[(b)] Suppose that {\rm {\bf (i)}} holds, with
$[a,b]=[0,\|x\|]$ for some  $x\in Z_*$,  and suppose that $x$
satisfies the further conditions:
\begin{description}
\item[(ii)] $\pair{x}{F_B}\subset B$ for each interval $B\in {\cal B}$;
\item[(iii)] $S_{F_B}^*x=x$ for $B\in {\cal B}$;
\item[(iv)] $\pair{x}{F_{(0,\|x\|]}^\perp}=0$.
\end{description}
Then $x=\int t\, dv_B(t)$.
\end{description}
\end{lemma}

\begin{proof} For the proof of (a), it suffices to show that for every
$\epsilon>0$, there is a $\delta>0$, such that
\begin{equation}\label{eq:Riemann}
\|S(P,T,f)-S(P',T',f)\|<\epsilon \mbox{ if }|P|,|P'|<\delta.
\end{equation}
By the uniform continuity of $f$, let $\delta>0$ correspondence to
a tolerance of $\epsilon/2$. If $|P|,|P'|<\delta$, then
$S(P,T,f)-S(P\cup P',T'',f)$, where $T''$ is any selection of
points, is of the form $\sum_1^m\alpha_jv_j$, where
$|\alpha_j|<\epsilon/2$ and $v_1,\ldots,v_m$ are orthogonal
geometric tripotents. Thus
\[
\|S(P,T,f)-S(P\cup P',T'',f)\|=\max_j |\alpha_j|<\epsilon/2
\]
and (\ref{eq:Riemann}) follows.

For the proof of (b), it suffices to prove that $x$ is the
weak*-limit of the Riemann sums corresponding to $f_0(t):=t$, for
by (a), $x$ will also be the norm limit. In what follows,
$F_{(0,\|x\|]}$ will be denoted by $F$. By (iii) and (iv)
\[
\pair{x}{Z_1(F)+Z_0(F)}=0.
\]
Also, each Riemann sum $\sum t_j v_{(s_{j-1},s_j]}\in U_2(F)$, so
\[
\pair{\sum t_j v_{(s_{j-1},s_j]}}{Z_1(F)+Z_0(F)}=0.
\]
Since $Z_2(F)=\mbox{sp}_\CC F$, it suffices to prove that for
every $\psi\in F$,
\[
\pair{x-S(P,T,f_0)}{\psi}\rightarrow 0\mbox{ as }|P|\rightarrow 0.
\]
Since $v_F=\sum v_i$ where $v_i=v_{(s_{i-1},s_i]}$, if $\psi\in
F\subset \oplus_i Z_2(v_i)\oplus\oplus_{i\ne j}[Z_1(v_i)\cap
Z_1(v_j)]$, then
\begin{eqnarray*}
1&=&\pair{v_F}{\psi}=\pair{v_F}{\sum P_2(v_i)\psi+\sum_{i\ne
j}P_1(v_i)P_1(v_j)\psi}\\
&=&\sum\pair{v_i}{P_2(v_i)\psi}\le \sum \|P_2(v_i)\psi\|=\|\sum
P_2(v_i)\psi\|\le \|\psi\|=1.
\end{eqnarray*}
Therefore
\begin{eqnarray*}
\psi&=&\sum\|P_2(v_i)\|\frac{P_2(v_i)}{\|P_2(v_i)\|}+\sum_{i\ne
j}P_1(v_i)P_1(v_j)\psi\\
&\in& \mbox{co}\, (F_{v_1}\cup\cdots\cup F_{v_n})+\oplus_{i\ne
j}[Z_1(v_i)\cap Z_1(v_j)].
\end{eqnarray*}
By (iii), $ \pair{x}{Z_1(F_B)}=0\mbox{ for every }B\in {\cal B}. $
Therefore $\pair{x}{\psi}=\pair{x}{\sum \lambda_i\psi_i}$, where
$\psi_i\in F_{v_i}, \lambda_i\ge 0,\sum \lambda_i=1$. Also,
$\pair{S(P,T,f_0)}{\psi}=\pair{\sum t_iv_i}{\sum \lambda_j
\psi_j}=\sum t_i\lambda_i$.

By (ii), $\pair{x}{\psi_i}\in (s_{i-1},s_i]$, so
\[
|\pair{x-S(P,T,f_0)}{\psi}|=|\sum \lambda_i
(\pair{x}{\psi_i}-t_i)|\le |P|.
\]
The lemma is proved.\qed\end{proof}

\medskip

Let us observe that if $Z$ is the dual of a \jbst\ $A$, then each
element $x\in A$ satisfies the conditions (i)-(iv) of
Lemma~\ref{lem:spectral}. Indeed, if $C$ denotes the JB*-subtriple
of $A$ generated by $x$, then $C$ is isometric to a commutative
C*-algebra and consists of norm limits of elements $p(x)$ where p
is an odd polynomial on $(0,\|x\|]$, cf. \cite[1.15]{Kaup83} and
\cite[p. 438]{AraKau01}; and if $W$ denotes the \jbwst\ generated
by $x$ in $A^{**}$, then $W$ is a commutative von Neumann algebra.
Thus, if $x=w|x|$ is the polar decomposition of $x$ in $W$, and
$|x|=\int\lambda\, de_\lambda$ is the spectral decomposition of
$|x|$ in $W$, and the face $F_B$ is defined as the face exposed by
the tripotent $we(B)\in A^{**}$, then the family $\{F_B:B\in {\cal
B}\}$ satisfies (i), as shown in \cite[Theorem 3.2]{Harris81}. It
also follows from \cite[Theorem 3.2]{Harris81} that for every
$\epsilon>0$, there is a partition of $[0,\|x\|]$ such that
$\|x-\sum t_j v_{(s_{j-1},s_j]}\|<\epsilon$. If $B$ is a
subinterval of $[0,\|x\|]$, and $\rho\in F_B$, then with
$v_j=v_{(s_{j-1},s_j]}$, $B_j=B\cap (s_{j-1},s_j]$, there exist
$\rho_k\in F_{B_k}$ (if $B_k\ne\emptyset$) and $\lambda_k\ge 0$
with $\sum \lambda_k=1$ such that $\pair{x}{\rho}$ is approximated
by
\[
\pair{\sum t_j
v_j}{\sum_{B_k\ne\emptyset}\lambda_k\rho_k}=\sum_{B_j\ne\emptyset}
t_j\lambda_j\in \mbox{co}\, (\cup_{B_j\ne\emptyset}B_j),
\]
proving (ii). Again, using \cite[Theorem 3.2]{Harris81} we shall
show that (iii) and (iv) hold.  Since $x$ is approximated in norm
by $\sum t_jv_j$, where $v_j=v_{(s_{j-1},s_j]}$, to prove (iii),
it suffices to prove that
$v_Bv_B^*v_jv_B^*v_B=v_Bv_B^*v_j=v_jv_B^*v_B$. Since $v_B=\sum
v_{B_j}$ where $B_j=B\cap (s_{j-1},s_j]$, it is trivial to check
that each of the terms $v_Bv_B^*v_jv_B^*v_B,\ v_Bv_B^*v_j,\
v_jv_B^*v_B$ collapses to $v_{B_j}$. Since the support of the
spectral measure of $|x|$ lies in $[0,\|x\|]$, (iv) also holds.

There is another property of elements of a \jbst\ that we need to
incorporate into our definition. It is based on the following
observation. If $x$ is an element of a \jbst\ $A$, let $f(x)$
denote the element of $C$ which is the norm limit of odd
polynomials $p_n$ which converge uniformly to $f\in
C_0([0,\|x\|)$, and let $\tilde{f}(x)=\int f(\lambda)\,
de_\lambda$. Since $p_n(x)=\tilde{p}_n(x)$,
\begin{eqnarray*}
f(x)-\tilde{f}(x)&=&f(x)-p_n(x)+\tilde{p}_n(x)-\sum
p_n(t_k)v_k\\
&+&\sum p_n(t_k)v_k-\sum f(t_k)v_k +\sum f(t_k)v_k -\tilde{f}(x),
\end{eqnarray*}
which shows that $\tilde{f}(x)=f(x)\in A$.

\begin{definition}\label{ssdef}{\rm
A strongly facially symmetric space $Z$ with a predual $Z_*$ is
{\it  strongly spectral} if, for every element $x \in Z_{\ast}$,
there exists a family $\{ F_{B}:B\in {\cal B} \}$ of norm closed
faces of the closed unit ball $Z_{1}$, where $\cal B$ is the set
of nonempty Borel subsets of $(0,\|x\|]$, satisfying (i)-(iv) in
Lemma~\ref{lem:spectral} and which also satisfies
\begin{description}
\item[(v)] For every $f\in C_0(0,\|x\|)$, the element $\int f\,
dv_B$ is weak*-continuous, that is, lies in $Z_*$.
\end{description}}
\end{definition}

Although somewhat complicated, this condition is precisely the
analogue of a strongly spectral compact base $K$ of a base normed
space $V$ given by Alfsen and Shultz in \cite{AS}.  There it is
given simply as the condition that in the order unit space
$V_{\ast}$ each element $a$ decomposes as an orthogonal difference
$a_{+}-a_{-}$ of two positive elements. Here orthogonal means that
$a_{+}$ and $a_{-}$ are supported on real spans of orthogonal
faces of $K$. Since $V_{\ast}$ is unital, the unit may be used
together with $a$ and this property to carve out an orthogonal
collection of faces similar to the one above, and a lattice of
orthogonal elements of $V_{\ast}$ which generate a space which is
isometric to a full space of continuous functions, and hence
closed under the continuous functional calculus. Since there is no
unit in our space $Z_{\ast}$, we must assume that elements $x \in
Z_{\ast}$ may be decomposed in the above fashion, and that the
resulting continuous functional calculus operates in $Z_*$. Note
that this is entirely a linear property, and has obvious quantum
mechanical significance. The faces $F_{B}$ are the states
corresponding to observations of some value in $B$ for the
observable $x$. The probability if this happening for a state
$\psi$ is $|\psi(v_{B})|$.

\medskip

We now have the following characterizations of JB*-triples. In
this characterization, the property JP must hold for all
orthogonal faces, not just extreme points. Thus it simply says
that the (necessarily commutative) product of the symmetries
$S_{F}$ and $S_{G}$ corresponding to orthogonal faces $F$ and $G$
is $S_{F \vee G}$.

\begin{theorem}\label{face} A Banach space $X$ is
isometric to a JB*-triple if and only if $X^*$ is an L-embedded,
locally base normed, strongly spectral, strongly facially
symmetric space which satisfies the pure state properties and JP.
\end{theorem}

 Before proving this theorem, we require one more lemma.

\begin{lemma}\label{lem:1208}
Let $Z:=X^*$ and $\Psi$ {\rm (}resp. $\Psi^\perp${\rm )} denote
the projection of $Z$ onto its atomic part $Z_{a}$ {\rm (}resp.
nonatomic part $Z_n${\rm )} given by Theorem \ref{atdec}. For any
norm exposed face $G\subset Z_1$, $G_a:=\Psi(G)\cap \partial Z_1$
and $G_n:=\Psi^\perp(G)\cap
\partial Z_1$ are faces in $Z_a$ and $Z_n$ respectively, and
\begin{equation}\label{eq:1208}
G=\mbox{co}\, (G_a\cup G_n).
\end{equation}
 Moreover, writing $G=F_w$ for some geometric tripotent $w$,
then $\Psi^*w$ is a geometric tripotent, and
\begin{equation}\label{eq:1209}
F_{\Psi^*w}=G_a.
\end{equation}
\end{lemma}
\begin{proof} To show that $G_a$ is a face in $(Z_a)_1$, let
$\lambda\rho+(1-\lambda)\sigma\in G_a$ where $\rho,\sigma\in
(Z_a)_1$. Then $\lambda\rho+(1-\lambda)\sigma=\Psi f$ for some
$f\in G$, and $f=\lambda\rho+(1-\lambda)\sigma+f_n$. Since
$\|f\|=1=\|\lambda\rho+(1-\lambda)\sigma\|$, $f_n=0$ and
$\rho,\sigma\in G$, $\|\rho\|=1=\|\sigma\|$, and $\rho\in G\cap
Z_a$, proving that $G_a$ is a face. Similarly for $G_n$.

If $f\in G$ has decomposition
$f=f_a+f_n=\|f_a\|\frac{f_a}{\|f_a\|}+\|f_n\|\frac{f_n}{\|f_n\|}$,
then since $G$ is a face,
$\frac{f_a}{\|f_a\|},\frac{f_n}{\|f_n\|}\in G$. This proves
$\subset$ in (\ref{eq:1208}). If $g_a:=\Psi g\in \Psi G\cap
\partial Z_1$ for some $g\in G$, then $\|g_a\|=1$ so $g=g_a=\Psi
g\in G$. A similar argument for $\Psi^\perp(G)\cap \partial Z_1$
proves $\supset$ in (\ref{eq:1208}).

To prove (\ref{eq:1209}), let $g\in \Psi(G)\cap\partial Z_1$. Then
$\pair{g}{\Psi^*w}=\pair{g}{w}=1$ so that $g\in F_{\Psi^*w}$.  On
the other hand, if $g\in F_{\Psi^*w}$, then
$1=\|g\|=\pair{g}{\Psi^*w}=\pair{\Psi g}{w}$ so that $\Psi g\in
F_w$. Since $g=\Psi g+\Psi^\perp g$ and $\|g\|=\|\Psi g\|$,
$\Psi^\perp g=0$, $\Psi g=g$ and $g\in \Psi(G)\cap \partial Z_1$.

It remains to show that $\Psi^*w$ is a geometric tripotent, that
is, $$\pair{\Psi^*w}{(G_a)^\perp}=0.$$ Note first that
$G^\perp=G_a^\perp\cap G_n^\perp$ by (\ref{eq:1208}). If $\rho\in
G_a^\perp$, $\pair{\Psi^*w}{\rho}=\pair{w}{\Psi(\rho)}$ and this
will be zero if $\Psi(\rho)\in G^\perp$. To prove this, first let
$\sigma\in G_a$. Then $\rho\perp\sigma$, hence
$\Psi(\rho)\perp\Psi(\sigma)$ and since $\Psi(\sigma)=\sigma$,
$\Psi(\rho)\in G_a^\perp$. Then $\Psi(\rho)\in G_a^\perp\cap
G_n^\perp=G^\perp$ as required. \qed\end{proof}

\medskip

\noindent{\it Proof of Theorem \ref{face}.}
 Assume that $Z=X^*$ is a strongly facially symmetric space satisfying the
hypotheses of the theorem.  Suppose $x$ is an element of
$X=Z_{\ast}$. By the spectral axiom and Lemma~\ref{lem:spectral},
there is an element $y\in X$ such that for $\epsilon>0$ there
exists $\delta>0$ such that, with $f_0(t)=t$ and $f_1(t)=t^3$,
\[
\|x-S(P,T,f_0)\|<\epsilon\mbox{ and }\|y-S(P',T',f_1)\|<\epsilon
\]
for all partitions $P,P'$ with mesh less than $\delta$. Fix a
common partition $P=\{s_0,\ldots,s_n\}$ with $|P|<\delta$, and
write $v_i=v_{(s_{i-1},s_i]}$ and $(v_i)_a=\Psi^*(v_i)$. Then by
(\ref{eq:1209}),
\[
\|\Psi^*(\hat x)-\sum t_i(v_i)_a\|<\epsilon\mbox{ and
}\|\Psi^*(\hat y)-\sum t_i^3(v_i)_a\|<\epsilon.
\]

Since in a \jbst, $\|\{aaa\}-\{bbb\}\|\le
\|a-b\|(\|a\|^2+\|a\|\|b\|+\|b\|^2)$, and since the $(v_i)_a$ are
orthogonal tripotents in the \jbwst\  $(Z_a)^*$, we have
\[
\|\{\Psi^*(\hat x),\Psi^*(\hat x),\Psi^*(\hat x)\}-\sum
t_i^3(v_i)_a\|<3\epsilon\|x\|^2,
\]
and therefore $\|\{\Psi^*(\hat x),\Psi^*(\hat x),\Psi^*(\hat
x)\}-\Psi^*(y)\|<\epsilon(3\|x\|^2+1)$. It follows that
$\Psi(\widehat{X})$  is a norm closed subspace of the JBW*-triple
$(Z_{a})^{\ast}$ that is closed under the cubing operation.  Hence
$\Psi(X)$ is a subtriple of $(Z_{a})^{\ast}$ as required.

The converse, that the dual $Z$ of a JB*-triple is a strongly
facially symmetric space satisfying the conditions of the theorem,
has already been mentioned above. That the spectral axiom is
satisfied was shown preceding Definition~\ref{ssdef}. The proofs
that it is a strongly facially symmetric locally base normed space can be
found in \cite{Pac}, the proofs that it satisfies the pure state
properties can be found in \cite{FriRus85bis}, the proof of the
L-embeddedness can be found in \cite{BarTim86}, and the proof of
FE can be found in \cite{ER}.\qed

\medskip

We can restate Theorem~\ref{face} from another viewpoint as
follows: for a Banach space $X$, its open unit ball is a bounded
symmetric domain if and only if $X^*$ is an L-embedded, locally base
normed, strongly spectral, neutral strongly facially symmetric
space which satisfies the pure state properties and JP.


\subsection{One-sided ideals in \csa s}

Proposition~\ref{April} and Theorem~\ref{C} below, together with
Theorem~\ref{theorem:5.5},  give facial and linear operator space
characterizations of C*-algebras and left ideals of C*-algebras.
This work was inspired by \cite{Blecherpp}, in which
Theorem~\ref{theorem:5.5} is used to characterize left ideals as
TRO's which are simultaneously abstract operator algebras with
right contractive approximate unit.

We start by motivating the main result of this subsection. Recall
that a TRO is made into a \jbst\ by symmetrizing the ternary
product.
\begin{remark}
If $J$ is a closed left ideal in a C*-algebra and $J$ possesses a
right identity $e$ of norm 1, then $J$ is a TRO and
$E:=\left[\begin{array}{c}0\\ e \end{array}\right]$ is a maximal
partial isometry in $M_{2,1}(J)$, that is,  $P_0(E)=0$.
\end{remark}
\begin{proof} By a remark of Blecher (see \cite[Lemma 2.9]{Blecherpp}),
$xe^*=x$ for all $x\in J$, so that $x=xe^*e$ and in particular,
$e$ is a partial isometry, and so is $E$.

For $\left[\begin{array}{c}x\\ y
\end{array}\right]\in M_{2,1}(J)$,
\begin{eqnarray*}
P_0(E) \left[\begin{array}{c}x\\ y \end{array}\right]&=&(I-EE^*)
\left[\begin{array}{c}x\\ y \end{array}\right](I-E^*E)\\
&=&\left[\begin{array}{cc}1 & 0\\ 0&1-ee^* \end{array}\right]
\left[\begin{array}{c}x \\ y\end{array}\right](I-e^*e)\\
&=&\left[\begin{array}{c}x(1-e^*e)\\ (1-ee^*)y(1-e^*e)
 \end{array}\right]=0. \qed
\end{eqnarray*}
\end{proof}

 Conversely, we have the following.

\begin{proposition}
Let $A$ be a TRO. Suppose there is a norm one element $x$ in $A$
such that the face in $(M_{2,1}(A)^*)_1$ exposed by
$$X:=\left[\begin{array}{c}0\\ x \end{array}\right]\in M_{2,1}(A)$$
 is maximal. Then $A$ is completely isometric
to a left ideal in a C*-algebra, which ideal contains a right
identity element.
\end{proposition}
\begin{proof}
 Let $B=M_{2,1}(A)$. If $V$ is the partial isometry in
$B^{**}$ such that $F_X=F_V$, then $X=V+P_0(V)^*X=V$, so that $x$
is a partial isometry in $A$, which we denote by $v$.

We next prove that $v$ is a right unitary in $A$; that is,
$x=xv^*v$, for all $x\in A$. Indeed, for $x\in A$,
\begin{eqnarray*}
D(V)\left[\begin{array}{c}x\\ 0
\end{array}\right]&=&\frac{1}{2}\left( \left[\begin{array}{c}0\\ v
\end{array}\right] \left[\begin{array}{c}0\\ v
\end{array}\right]^* \left[\begin{array}{c}x\\ 0
\end{array}\right] + \left[\begin{array}{c}x\\ 0
\end{array}\right] \left[\begin{array}{c}0\\ v
\end{array}\right]^*
\left[\begin{array}{c}0\\ v \end{array}\right]\right)\\
&=&\left[\begin{array}{c}xv^*v/2\\ 0 \end{array}\right],
\end{eqnarray*}
and
\[
P_2(V)\left[\begin{array}{c}x\\ 0 \end{array}\right]=
\left[\begin{array}{c}0\\ v \end{array}\right]
\left[\begin{array}{c}0\\ v \end{array}\right]^*
\left[\begin{array}{c}x\\ 0 \end{array}\right]
\left[\begin{array}{c}0\\ v \end{array}\right]^*
\left[\begin{array}{c}0\\ v \end{array}\right]=0.
\]
Since $P_1(V)\left[\begin{array}{c}x\\ 0 \end{array}\right]=
\left[\begin{array}{c}x\\ 0 \end{array}\right]$, and
\[
\left[\begin{array}{c}xv^*v/2\\ 0 \end{array}\right]=
D(V)\left[\begin{array}{c}x\\ 0 \end{array}\right]=
P_2(V)\left[\begin{array}{c}x\\ 0 \end{array}\right]+\frac{1}{2}
P_1(V)\left[\begin{array}{c}x\\ 0 \end{array}\right]=\frac{1}{2}
\left[\begin{array}{c}x\\ 0 \end{array}\right].
\]

We next show that the map $\psi:a\mapsto av^*$ is a complete
isometry of $A$ onto a closed left ideal $J$ of the C*-algebra
$\overline{AA^*}$ and $vv^*$ is a right identity of $J$. In the
first place, since
$\|\psi(x)\|^2=\|xv^*\|^2=\|(xv^*)(xv^*)^*\|=\|xv^*vx^*\|=\|xx^*\|$,
$\psi$ is an isometry. By the same argument, with $W=\mbox{diag}\,
(v,v,\ldots,v)$, for $X\in M_n(A)$, $\|XW^*\|=\|X\|$, so that
$\psi$ is a complete isometry.

If $c\in \overline{AA^*}$ is of the form $c=ab^*$ with $a,b\in A$,
and $y\in J:=\psi(A)$, say $y=xv^*$, then $cy=ab^*xv^*\in Av^*=J$.
By taking finite sums and then limits, $J$ is a left ideal in $C$.
Finally, with $e=vv^*$ and $y=xv^*\in J$,
$ye=xv^*vv^*=xv^*=y$.\qed\end{proof}

\medskip


For the general case we have the following result.
\begin{theorem}\label{C}
Let $A$ be a TRO. Then $A$ is completely isometric to a left ideal
in a C*-algebra if and only if there exists a  convex set
$C=\{x_\lambda:\lambda\in \Lambda\}\subset A_1$  such that the
collection of faces
$$F_\lambda:=F_{\left[\begin{array}{c}0\\
x_\lambda/\|x_\lambda\|
\end{array}\right]}\subset M_{2,1}(A)^*,$$ form a directed set with
respect to containment, $F:=\sup_{\lambda}F_\lambda$ exists, and
\begin{description}
\item[(a)] The set $\{\left[\begin{array}{c}0\\ x_{\lambda}
\end{array}\right]:\lambda\in\Lambda\}$ separates the points of $F$;
\item[(b)] $F^\perp=0$ {\rm (}that is, the partial isometry $V\in
(M_{2,1}(A))^{**}$ with $F=F_V$ is maximal{\rm )};
\item[(c)] $\pair{F}{\left[\begin{array}{c}0\\ x_\lambda
\end{array}\right]}\ge 0$ for all $\lambda\in\Lambda$;
\item[(d)] $S_F^*\left(\left[\begin{array}{c}0\\ x_\lambda
\end{array}\right]\right)=\left[\begin{array}{c}0\\ x_\lambda
\end{array}\right]$ for all $\lambda\in\Lambda$.
\end{description}
\end{theorem}
\begin{proof}  We first assume that we have a closed left ideal $L$ in a
C*-algebra $B$. In this part of the proof, to avoid confusion with
dual spaces, we denote the involution in $B$ by $x^\sharp$. The
set of positive elements of the open unit ball of the C*-algebra
$L \cap L^{\sharp}$, which we will denote by
$(u_\lambda)_{\lambda\in\Lambda}$, is a contractive right
approximate unit for $L$.
 Let
$u=\mbox{w}^*$-$\lim u_\lambda\in B^{**}$. Identifying
$L^{\ast\ast}$ with $B^{\ast\ast}u$, we now verify the properties
(a)-(d).

For each $\lambda$,
$u_\lambda/\|u_\lambda\|=v_\lambda+v_\lambda^0$ where
$v_{\lambda}=\mbox{w}^*\mbox{-}\lim (u_\lambda/\|u_\lambda\|)^{n}$
is the support projection of $u_{\lambda}/\|u_\lambda\|$, that is,
$F_{u_{\lambda}/\|u_\lambda\|}=F_{v_{\lambda}}\subset B^*$, and
$v_\lambda^0$ is an element orthogonal to $v_\lambda$. Since
$u_\lambda \uparrow u$, $u=\sup_\lambda
r(u_\lambda/\|u_\lambda\|)$, where $r(u_\lambda/\|u_\lambda\|)$ is
the range projection of $u_\lambda/\|u_\lambda\|$. For each fixed
$\mu\in\Lambda$, we apply the functional calculus to
$u_\mu/\|u_\mu\|$ as follows. Let $f_n(0)=0,\ f_n(t)=1$ on
$[1/n,1]$ and linear on $[0,1/n]$. Then $f_n(u_\mu/\|u_\mu\|)\in
(L\cap L^\sharp)^+_1$ and so as above
$f_n(u_\mu/\|u_\mu\|)=v_{\lambda(\mu,n)}+v_{\lambda(\mu,n)}^0$ and
$\sup_n v_{\lambda(\mu,n)}=r(u_\mu/\|u_\mu\|)$. Therefore
\[
u=\sup_\mu r(u_\mu/\|u_\mu\|)=\sup_\mu\sup_n v_{\lambda(\mu,n)}\le
\sup_\lambda v_\lambda=v\mbox{ say}.
\]
On the other hand, since $v_\lambda\le
(1+\frac{1-\|u_\lambda\|}{\|u_\lambda\|})u$, it follows that $v\le
u$ and therefore $u=v$.

It is clear that $$F_\lambda=F_{\left[\begin{array}{c}0\\
u_\lambda/\|u_\lambda\|
\end{array}\right]}=
F_{\left[\begin{array}{c}0\\
v_\lambda
\end{array}\right]}
\subset F_{\left[\begin{array}{c}0\\ u
\end{array}\right]},$$
and therefore that $\sup_\lambda F_\lambda$ exists. We show that
it equals $F_{\left[\begin{array}{c}0\\ u
\end{array}\right]}$.
Suppose that for some $a,b\in B^{**}$, $F_\lambda\subset F_{\left[\begin{array}{c}a\\
b
\end{array}\right]}$ for every $\lambda$. This
is equivalent to
 $$\left[\begin{array}{c}0\\ v_\lambda
\end{array}\right]=Q(\left[\begin{array}{c}0\\ v_\lambda
\end{array}\right])\left[\begin{array}{c}a\\ b
\end{array}\right]=\left[\begin{array}{c}0\\v_\lambda
b^{\ast}v_\lambda
\end{array}\right],$$ or $v_\lambda b^*v_\lambda=v_\lambda$. On the
other hand, since
$v_\lambda^0=u_\lambda/\|u_\lambda\|-v_\lambda\rightarrow 0$,
\[
u_\lambda b^*
u_\lambda=\|u_\lambda\|^2(v_\lambda+v_\lambda^0)b^*(v_\lambda+v_\lambda^0)\rightarrow
u,
\]
so that $ub^*u=u$ and as above, $F_{\left[\begin{array}{c}0\\ u
\end{array}\right]}\subset F_{\left[\begin{array}{c}a\\ b
\end{array}\right]}$, proving that $\sup_\lambda F_\lambda=F_{\left[\begin{array}{c}0\\ u
\end{array}\right]}$.

\medskip

 Let us now prove (a). Since $u_\lambda\uparrow u$, the convergence is strong
convergence. We claim first that $B\cap B^{**}_2(u)$ is
weak*-dense in $B_2^{**}(u)$. Indeed with $x\in B_2^{**}(u)$ of
norm 1, there is a net $b_\alpha\in B$ with $b_\alpha\rightarrow
x$ strongly. Then $u_\lambda b_\alpha u_\mu\rightarrow uxu=x$, and
since $u_\mu u=u_\mu r(u_\mu)u=u_\mu  r(u_\mu)=u_\mu$, $u_\lambda
b_\alpha u_\mu\in B\cap B^{**}_2(u)$, proving the claim. We claim
next that $L\cap L^\sharp=B_2^{**}(u)\cap B$. If $y\in L\cap
L^\sharp$, then $y=bu=uc^\sharp$ for some $b,c\in B^{**}$, hence
$y=uyu\in B_2^{**}\cap B$. Since $B_2^{**}(u)=uB^{**}u\subset
B^{**}u=L^{**}$, we have $B_2^{**}(u)\cap B\subset L^{**}\cap
B=L$. If $x\in B_2^{**}(u)\cap B$, then $x^\sharp\in
B_2^{**}(u)\cap B$, proving that $x\in L\cap L^\sharp$.

Let $M$ denote the TRO $\left[\begin{array}{c}L\\
L\end{array}\right]$. Let $f,g$ be two elements of
$F_{\left[\begin{array}{c}0\\ u
\end{array}\right]}$ which are not separated by
$\left[\begin{array}{c}0\\ C
\end{array}\right]$.  It follows that $\left[\begin{array}{c}0\\ C
\end{array}\right]$ annihilates $f-g \in
M^{\ast}_{2}(\left[\begin{array}{c}0\\ u
\end{array}\right])=\mbox{sp}_\CC F_{\left[\begin{array}{c}0\\ u
\end{array}\right]}$.
This contradicts the fact, implicit in the preceding paragraph,
that the linear span of $C$ is w*-dense in $
L^{\ast\ast}_{2}(u)=B_2^{**}(u)$. This proves (a).

\medskip

To prove (b), it suffices to show that for $a,b\in L^{**}$,
\[
\left(1-\left[\begin{array}{c}0\\
u\end{array}\right]\left[\begin{array}{c}0\\
u\end{array}\right]^*\right)
\left[\begin{array}{c}a\\
b\end{array}\right]
\left(1-\left[\begin{array}{c}0\\
u\end{array}\right]^*\left[\begin{array}{c}0\\
u\end{array}\right]\right)=0.
\]

This reduces to
\[
\left[\begin{array}{cc} 1&0\\
0&1-u
\end{array}\right]
\left[\begin{array}{c}a\\
b\end{array}\right](1-u)=\left[\begin{array}{c}a(1-u)\\
(1-u)b(1-u)\end{array}\right]=0,
\]
which is true since $u$ is a right identity for $L^{**}$.

To prove (c), let $N$ denote the TRO $\left[\begin{array}{c}B\\
B\end{array}\right]$. Note that $F_{\left[\begin{array}{c}0\\
u\end{array}\right]}$ is the normal state space of the von Neumann
algebra
$N^{**}_2\left(\left[\begin{array}{c}0\\
u\end{array}\right]\right)$ and that
\[
\left[\begin{array}{c}0\\
u_\lambda
\end{array}\right]=
\left[\begin{array}{c}0\\
\sqrt{u_\lambda}
\end{array}\right]
[0, u]
\left[\begin{array}{c}0\\
\sqrt{u_\lambda}
\end{array}\right]
\]
is the square of the self-adjoint element
\[
\left[\begin{array}{c}0\\
\sqrt{u_\lambda}
\end{array}\right]^\sharp=
\left[\begin{array}{c}0\\
u\end{array}\right] [0, \sqrt{u_\lambda}]
\left[\begin{array}{c}0\\
u\end{array}\right]=
\left[\begin{array}{cc} 0&0\\
0&\sqrt{u_\lambda}
\end{array}\right]
\left[\begin{array}{c}0\\
u\end{array}\right]=\left[\begin{array}{c}0\\
\sqrt{u_\lambda}\end{array}\right].
\]
Hence (c) follows.

From the proof of (c), $\left[\begin{array}{c}0\\
u_\lambda
\end{array}\right]\in N^{**}_2\left(\left[\begin{array}{c}0\\
u\end{array}\right]\right)$, so it is fixed by $S_F^*$.

\medskip

 To prove the converse, assume that $A$ is a TRO satisfying
the conditions of the theorem. Let $B$ denote the TRO
$\left[\begin{array}{c}A\\A\end{array}\right]$. As in the first
part of the proof, for each $\lambda$, there exists a partial
isometry $v_{\lambda}\in A^{\ast\ast}$ and an element
$v^{0}_{\lambda} \in A_{0}^{\ast\ast}(v_{\lambda})$ such that
$\left[\begin{array}{c}o\\x_{\lambda}/\|x_\lambda\|\end{array}\right]=\left[\begin{array}{c}0\\v_{\lambda}\end{array}\right]+\left[\begin{array}{c}0\\v^{0}_{\lambda}\end{array}\right]$
and
$F_{\left[\begin{array}{c}0\\x_{\lambda}/\|x_\lambda\|\end{array}\right]}=F_{\left[\begin{array}{c}0\\v_{\lambda}\end{array}\right]}$.
Since $\sup_\lambda
F_{\left[\begin{array}{c}0\\v_{\lambda}\end{array}\right]}=F$
exists, let $F=F_{\left[\begin{array}{c}u\\v\end{array}\right]}$
with $\left[\begin{array}{c}u\\v\end{array}\right]$ a partial
isometry in $(M_{2,1}(A))^{**}$. We shall show that $u=0$ and
hence that $v$ is a partial isometry. In the first place,
$P_2(\left[\begin{array}{c}0\\v_\lambda\end{array}\right])(\left[\begin{array}{c}u\\v\end{array}\right])
=\left[\begin{array}{c}0\\v_\lambda\end{array}\right]$, which
reduces to $P_2(v_\lambda)v=v_\lambda$. Since
$\left[\begin{array}{c}0\\v\end{array}\right]$ is the image of
$\left[\begin{array}{c}u\\v\end{array}\right]$ under a contractive
projection, $\|v\|\le 1$, and therefore $P_1(v_\lambda)^*v=0$ (by
\cite[Lemma 1.5]{FriRus85bis}). Thus $v=v_\lambda+v_\lambda^0$
with $v_\lambda^0$ orthogonal to $v_\lambda$, and it follows that
the support partial isometry $u(v)$ of the element $v\in A^{**}$
satisfies $u(v)\ge v_\lambda$.  It follows that
$\left[\begin{array}{c}0\\v_\lambda\end{array}\right] \le
\left[\begin{array}{c}0\\u(v)\end{array}\right]$ and since
$\left[\begin{array}{c}u\\v\end{array}\right]$ is the least upper
bound, we have $\left[\begin{array}{c}u\\v\end{array}\right] \le
\left[\begin{array}{c}0\\u(v)\end{array}\right]$. Thus
\[
\left[\begin{array}{c}u\\v\end{array}\right] = P_2(
\left[\begin{array}{c}0\\u(v)\end{array}\right] )
\left[\begin{array}{c}u\\v\end{array}\right]
=\left[\begin{array}{c}0\\P_2(u(v))v\end{array}\right] =
\left[\begin{array}{c}0\\v\end{array}\right],
\]
showing that $u=0$.

Conditions (b) and (d) imply that
$\left[\begin{array}{c}0\\x_{\lambda}\end{array}\right]$ lies in
the von Neumann algebra
$B_{2}^{\ast\ast}(\left[\begin{array}{c}0\\v\end{array}\right])$
while condition (c) implies that
$\left[\begin{array}{c}0\\x_{\lambda}\end{array}\right] \geq 0$ in
that von Neumann algebra. In particular,
$\left[\begin{array}{c}0\\x_\lambda\end{array}\right]$ is
self-adjoint; $vx_\lambda^*v=x_\lambda$. We claim that condition
(a) implies that $\left[\begin{array}{c}0\\C\end{array}\right]$
cannot annihilate any non-zero element of $B^{\ast}_{2}(
\left[\begin{array}{c}0\\v\end{array}\right])$. Indeed, suppose
$\left[\begin{array}{c}0\\C\end{array}\right](\psi_{1}-\psi_{2})=0$
where $\psi_{1}-\psi_{2}$ is the Jordan decomposition of a
functional $\psi$ in the self adjoint part of
$B^{\ast}_{2}(\left[\begin{array}{c}0\\v\end{array}\right])$. Note
that since $\{v_\lambda\}$ is directed, and $v_\lambda\le
x_\lambda\le v$, it follows that
$\|\psi_{1}\|=\left[\begin{array}{c}0\\v\end{array}\right](\psi_{1})=\sup
\left[\begin{array}{c}0\\C\end{array}\right](\psi_1)= \sup
\left[\begin{array}{c}0\\C\end{array}\right](\psi_2)=
\left[\begin{array}{c}0\\v\end{array}\right](\psi_{2})=\|\psi_{2}\|$
and this contradicts (a), as $\psi_1/\alpha,\psi_2/\alpha\in F$,
where $\alpha$ is the common norm of $\psi_1$ and $\psi_2$. It
follows that the bipolar
$(\left[\begin{array}{c}0\\C\end{array}\right]_0)^0=B_{2}^{\ast\ast}
( \left[\begin{array}{c}0\\v\end{array}\right])$. Consequently,
the w* closure of $\mbox{sp}_{\bf C}C$ is $A^{\ast\ast}_{2}(v)$
and since the norm closure of a convex set is the same as its weak
closure
\[
A\cap A_2^{**}(v)=A\cap \overline{\mbox{sp}\, C}^{\mbox{w*}}=A\cap
\overline{\mbox{sp}\, C}^{\|\cdot\|}=\overline{\mbox{sp}\,
C}^{\|\cdot\|}
\]
is a C*-subalgebra of $A^{\ast\ast}_{2}(v)$.

We are now in a position to show that $A$ is completely isometric
to a left ideal of a C*-algebra. Exactly as in the proof of the
right unital case we have $A\subset Av^{\ast}v$.  We define a map
$\Psi:A \rightarrow AA^{\ast}$ by $\Psi(a)=av^{\ast}$. The crux of
the matter is to show that the range of $\Psi$ lies in
$AA^{\ast}$. If that is the case, then since for $X,Y,Z\in
M_n(A)$, with $D=\mbox{diag}(v^*,\ldots,v^*)$,
\[
XY^*ZD=XD(YD)^* ZD,
\]
$\psi$ is a complete isometry. Moreover, if $b,c\in A$ then
$(bc^*)av*=(bc^*a)v^*$ shows that the range of $\psi$ is a left
ideal. It remains to show that $Av^*\subset AA^*$.

Note first that,  for $a\in A$, $av^*x_\lambda\in A$, since
\begin{eqnarray*}
av^*x_\lambda&=&av^*x_\lambda^{1/2}\cdot
x_\lambda^{1/2}=av^*(v(x_\lambda^{1/2})^*v)v^*x_\lambda^{1/2}\\
&=&
(av^*v)(x_\lambda^{1/2})^*(vv^*x_\lambda^{1/2})=a(x_\lambda^{1/2})^*x_\lambda^{1/2}\in
A.
\end{eqnarray*}
Next, since $v$ belongs to the w*-closure of $ \mbox{sp}_\RR\, C$,
and for each $a\in A$, $\{av^*y:y\in  \mbox{sp}_\RR\, C\}$ is a
convex subset of $A$ (since $av^*y=\sum\alpha_i av^*x_{\lambda_i}
=\sum\alpha_i a(x_{\lambda_i}^{1/2})^*x_{\lambda_i}^{1/2}\in A$),
it follows that $a$ belongs to the norm closure of  $\{av^*y:y\in
\mbox{sp}_\RR\, C\}$.  Now
$va^*av^*y=va^*av^*vy^*v=va^*ay^*v=(ya^*a)^\sharp\in A\cap
A_2^{**}(v)$ and therefore $va^*a$ belongs to the norm closure of
the set $\{va^*av^*y:y\in  \mbox{sp}_\RR\, C\}$ and hence
$va^*a\in A$.  Using the triple functional calculus in the TRO $A$
(see \cite{NeaRuspj}) , we have
$$av^*=a^{1/3}(a^{1/3})^*a^{1/3}v^*=a^{1/3}(v(a^{1/3})^*a^{1/3})^*\in
AA^*.\qed $$ \end{proof}

 \medskip

 In Theorem~\ref{C}, the elements $x_{\lambda}$  represent a
right approximate unit cast in purely linear terms.  Similar
language can be used to characterize C*-algebras.

\begin{proposition}\label{April}
Let $A$ be a TRO. Then $A$ is completely isometric to a unital
C*-algebra if and only if there is a norm one element $x$ in $A$
such that the complex linear span $\mbox{sp}\, _{\bf C}(F))$ of
the face $F$ in $A^{\ast}$ exposed by $x$ coincides with
$A^{\ast}$.
\end{proposition}

Note that a characterization of non-unital C*-algebras can also be
given with obvious modifications as in Theorem~\ref{C}.

From another viewpoint, we have characterized TRO's $A$ up to
complete isometry by {\it facial properties} of $M_{n}(A)^{\ast}$,
since by Theorem~\ref{theorem:5.5}, this is equivalent to finding
an isometric characterization of JB*-triples $U$ in terms of
facial properties of $U^{\ast}$. This is exactly what we have done
in Theorem~\ref{face}, which is the non-ordered version of
Alfsen-Shultz's facial characterization of state spaces of
JB-algebras  in the pioneering paper \cite{AS}.

\end{document}